# The Collatz problem $(a \cdot q+-1, a=1,3,5,\ldots)$ from the point of view of transformations of Jacobsthal numbers

Petro Kosobutskyy

Lviv Polytechnic National University Lviv Ukraine,

e-mail: petkosob@gmail.com

**Abstract:** In this work studies the general model $(a \cdot q+-1)$, $a=1,3,5,\ldots$ of the transformation of integer positive $q$ numbers in the directions of increasing (Jacobsthal transformation) and decreasing (Collatz conjecture) the degree natural $n$ numbers of the $Q \cdot 2^n$ $(Q=1,3,5,\ldots)$ sequences. It is shown that nodes are formed on the powers $Q \cdot 2^n$ under the condition $(Q \cdot 2^n +-1)/a=Integer$, from which the graphs are in the form of a tree (Jacobsthal tree). It is shown also that the trajectories of conjectures $(a \cdot q+-1)$ are formed on the Jacobsthal tree in the reverse direction of branches of $Q \cdot 2^n$ graphs. Reasoned analytical relations that allow you to calculate the parameters of the formation of periodic cycles

**PetroKosobutskyy** **ORCID iD** https://orcid.org/0000-0003-4319-7395

## Introduction

In mathematics, the Collatz problem [1] is known, in which a sequence of natural numbers $q \in \mathbb{N}$ is generated. Members of such a sequence are calculated according to the rule: let the number $q$ be even, then the next member of the sequence is equal to $q/2$, otherwise if $q$ is odd ($q_{odd}$), then the next member of the sequence is calculated as $C_{3q+1} = 3q+1$.The sequence for which this rule is fulfilled will be called the Collatz sequence, and its graph will be called the Collatz trajectory ($CT_{3q+1}$).

The analysis of studies [2-6] shows that the Collatz problem remains relevant, in which the transformation functions of positive integers have the form:

$$C_{3q\pm1} = \; if \quad q \equiv 0 \mod 2 \quad then \quad \frac{q}{2} \quad else \quad 3q \pm 1 \text{ , } q \in \mathbb{N}. \qquad (1)$$

However, the regularities of the conjecture $3q-1$ are fundamentally different from the classical Collatz conjecture $C_{3q+1}$. If in the direction decrease ($n \to 0$) of $n$-th power of two in $\theta \cdot 2^n$ the Collatz trajectories ($CT_{3q+1}$) lead to unity for all numbers $q \in [1,+\infty)$, then $CT_{3q-1}$ for $q$=5,7,17 , unity does not reach [7-8].

This work is a continuation of the work [9]. In this work, for the first time, based on the laws of transformation of Jacobsthal numbers, a model of the general Collatz problem

$$C_{a\cdot q\pm 1} = \ if \quad q \equiv 0 \mod 2 \quad then \quad \frac{q}{2} \quad else \quad a\cdot q \pm 1 \ , a = 1,3,5,... \in \mathbb{N}_{odd}, \ \mathbb{N} = \mathbb{N}_{odd} + \mathbb{N}_{even} \quad (2)$$

of the conjecture of natural numbers is developed. The first such studies were initiated by the author in [10]. Certain cases of conjectures (1) were studied in [3,5-8,11-17]. Jacobsthal numbers in Collatz conjectures were considered in [18]. Other mathematical and physical problems in the Collatz problem were addressed in [19-20]. The analysis of the results of fundamental studies of the Collatz problem [21-23], made in [2].

**Basic results, and their discussion**

In problem (2), the function $C_{a\cdot q_{odd}\pm 1} = a\cdot q_{odd} \pm 1$ conjecture the odd number $q_{odd}$ once into the even number $q_{even}$, which, by halving ($q/2$) $k_{q/2}$ times, transforms into the value of the first element $\theta\cdot 2^0$ of the binary sequence $\{\theta\cdot 2^n\}$. The number of iterations to transform an even number into an odd one does not depend on the form of the function $C_{a\cdot q\pm 1}$ i.e., on the value of the parameter $a$. Therefore, we construct the diagram of iterations of even numbers (highlighted in gray) of type

$$k_{q/2} = \underbrace{((q_{even}/2)/2)/2... \rightarrow q_{odd}}_{k\ iteration} \quad (3a)$$

| | | | | | | | | | | | | | | | | | | | | | | |
|---|---|---|---|---|---|---|---|---|---|---|---|---|---|---|---|---|---|---|---|---|---|---|
| | | | | | | | | | $k_{q/2}=1$ | | | | | | $k_{q/2}=2$ | | | | | | $k_{q/2}>2$ | |
| | 0 | 2 | 4 | 6 | 8 | | 0 | 2 | 4 | 6 | 8 | | 0 | 2 | 4 | 6 | 8 | | 0 | 2 | 4 | 6 | 8 |
| 00 | 0 | 1 | 2 | 1 | 3 | | | 1 | | 1 | | | | | 2 | | | | 0 | | | | 3 |
| 10 | 1 | 2 | 1 | 4 | 1 | | 1 | | 1 | | 1 | | | 2 | | | | | | | | 4 | |
| 20 | 2 | 1 | 3 | 1 | 2 | | | 1 | | 1 | | | 2 | | | | 2 | | | | 3 | | |
| 30 | 1 | 5 | 1 | 2 | 1 | | 1 | | 1 | | 1 | | | | | 2 | | | | 5 | | | |
| 40 | 3 | 1 | 2 | 1 | 4 | | | 1 | | 1 | | | | | 2 | | | | 3 | | | | 4 |
| 50 | 1 | 2 | 1 | 3 | 1 | | 1 | | 1 | | 1 | | | 2 | | | | | | | | 3 | |
| 60 | 2 | 1 | 6 | 1 | 2 | | | 1 | | 1 | | | 2 | | | | 2 | | | | 6 | | |
| 70 | 1 | 3 | 1 | 2 | 1 | | 1 | | 1 | | 1 | | | | | 2 | | | | 3 | | | |
| 80 | 4 | 1 | 2 | 1 | 3 | = | | 1 | | 1 | | + | | | 2 | | | + | 4 | | | | 3 |
| 90 | 1 | 2 | 1 | 5 | 1 | | 1 | | 1 | | 1 | | | 2 | | | | | | | | 5 | |
| 100 | 2 | 1 | 3 | 1 | 2 | | | 1 | | 1 | | | 2 | | | | 2 | | | | 3 | | |
| 110 | 1 | 4 | 1 | 2 | 1 | | 1 | | 1 | | 1 | | | | | 2 | | | | 4 | | | |
| 120 | 3 | 1 | 2 | 1 | 7 | | | 1 | | 1 | | | | | 2 | | | | 3 | | | | 7 |
| 130 | 1 | 2 | 1 | 3 | 1 | | 1 | | 1 | | 1 | | | 2 | | | | | | | | 3 | |
| 140 | 2 | 1 | 4 | 1 | 2 | | | 1 | | 1 | | | 2 | | | | 2 | | | | 4 | | |
| 150 | 1 | 3 | 1 | 2 | 1 | | 1 | | 1 | | 1 | | | | | 2 | | | | 3 | | | |
| 160 | 5 | 1 | 2 | 1 | 3 | | | 1 | | 1 | | | | | 2 | | | | 5 | | | | 3 |
| 170 | 1 | 2 | 1 | 4 | 1 | | 1 | | 1 | | 1 | | | 2 | | | | | | | | 4 | |
| 180 | 2 | 1 | 3 | 1 | 2 | | | 1 | | 1 | | | 2 | | | | 2 | | | | 3 | | |
| ... | | | | | | | | | | | | | | | | | | | | | | | |

(3b)

We see that in the $n \rightarrow 0$ direction, even numbers are divided into two groups:

$k_{q/2}$=1: 2→1, 10→5, 6→3, 14→7, 22→11, 30→15, …

$k_{q/2}$=2: 4→2→1, 12→6→3, 20→10→5, 28→14→7, 36→18→9,…

by directions $k = const$, and one direction

$k_{q/2}$=3,4,5,…: 8,16, 24,32,40, 48,56,64,72,80, 88,96,104,112,120,…

with variable value $k_{q/2}$. In addition, as follows from (4), two adjacent even numbers 62 and 64 have a fundamentally different number of iterations $k_{q/2}$ of transformation into an odd number: 62→62/2=31, i.e. $k_{q/2}$ =1, while for 64→ 32→16→8→4→ 2→1, i.e. for $k_{q/2}$ =6. This means that a comparative analysis of the laws of the Collatz transformation of two close numbers can lead to a incorrect conclusion.

Below is an expanded mapping of fragment diagram (3b) to odd numbers $q_{odd}$ :

Diagram of split iterations $k_{q/2} > 2$ of $q_{even} \to q_{odd}$ conjecture

| | 0 | 2 | 4 | 6 | 8 | | 0 | 2 | 4 | 6 | 8 | | 0 | 2 | 4 | 6 | 8 | | 0 | 2 | 4 | 6 | 8 |
|---|---|---|---|---|---|---|---|---|---|---|---|---|---|---|---|---|---|---|---|---|---|---|---|
| 00 | | 1 | 1 | 3 | 1 | | | 1 | | 3 | | | | | 1 | | | | 0 | | | | 1 |
| 10 | 5 | 3 | 7 | 1 | 9 | | 5 | | 7 | | 9 | | | 3 | | | | | | | | 1 | |
| 20 | 5 | 11 | 3 | 13 | 7 | | | 11 | | 13 | | | 5 | | | | 7 | | | | 3 | | |
| 30 | 15 | 1 | 17 | 9 | 19 | | 15 | | 17 | | 19 | | | | | 9 | | | | 1 | | | |
| 40 | 5 | 21 | 11 | 23 | 3 | | | 21 | | 23 | | | | | 11 | | | | 5 | | | | 3 |
| 50 | 25 | 13 | 27 | 7 | 29 | = | 25 | | 27 | | 29 | + | | 13 | | | | + | | | | 7 | |
| 60 | 15 | 31 | 1 | 33 | 17 | | | 31 | | 33 | | | 15 | | | | 17 | | | | 1 | | |
| 70 | 35 | 9 | 37 | 19 | 39 | | 35 | | 37 | | 39 | | | | | 19 | | | | 9 | | | |
| 80 | 5 | 41 | 21 | 43 | 11 | | | 41 | | 43 | | | | | 21 | | | | 5 | | | | 11 |
| 90 | 45 | 23 | 47 | 3 | 49 | | 45 | | 47 | | 49 | | | 23 | | | | | | | | 3 | |
| ... | | | | | | | | | | | | | | | | | | | | | | | |

We see that $k_{q/2}$=1 та $k_{q/2}$=2 conjecture of (3a) are displayed by a of type 1,3,5,7,9,11,13,15,17.19, … sequence. The structure of the constructed matrices does not depend on the type $a \cdot q_{odd} \pm 1$ conjecture of numbers $q_{odd}$ . Iteration $k_{q/2} > 2$ correspond to conjectures $\theta \cdot (((2^n/2)/2)/2...) \to \theta \cdot 2^0$.

The process (2) is realized in the $n \to 0$ direction, and in the reverse $n \to \infty$ direction, the tree of Jacobsthal branches is branched from the sequences $\theta \cdot 2^n$ nodes. Therefore, we justify the recurrent numbers of the nodes' in the binary sequence for specific models $a \in \mathbb{N}_{\text{odd}}$ :

$$a = 1: \ \{\theta \cdot 2^n\}, a = 3: \ \frac{1}{3}\{\theta \cdot 2^n\} = Integer, \quad a = 5: \ \frac{1}{5}\{\theta \cdot 2^n\} = Integer,$$
$$a = 7: \ \frac{1}{7}\{\theta \cdot 2^n\} = Integer, \quad a = 9: \ \frac{1}{9}\{\theta \cdot 2^n\} = Integer, \quad .... \qquad (4)$$

Let us show that parametrized $a, \theta$ integers of the type

$$\mathrm{K}^{\pm}_{a,\theta,n} = \frac{1}{a}\left[\theta \cdot 2^n \pm (-1)^n\right], \qquad (5)$$

form the so-called nodes of sequences (4), thanks to which in the direction $n \to \infty$ the even numbers of $q \in \mathbb{N}$ are transformed into odd ones and, conversely, in the reverse direction $n \to 0$ the odd numbers are conjectured into even ones.

According to the parameter $a$, the numbers (6) are structured with a period $T_a$ by powers of $n$ in sequence $\theta \cdot 2^n$. If $a = 1$, then $T_{a=1} = 2^0$. For $a = 3$, $T_{a=3} = 2^1$, the structuring of numbers (5) is analogous to the bisection by even and odd powers $n$. For $a = 5$, the period is equal to $T_{a=5} = 2^2$ and so on.

Let's divide the numbers (5) into two groups based on their signs:

$$\mathrm{K}^{-}_{a,\theta,n} = \frac{1}{a}\left[\theta \cdot 2^n - (-1)^n\right] \quad (a) \quad \text{and} \quad \mathrm{K}^{+}_{a,\theta,n} = \frac{1}{a}\left[\theta \cdot 2^n + (-1)^n\right], \quad (b) \quad (6)$$

after which we apply the bisection by powers of

$$n = N_{even} + N_{odd} = r + s \qquad (7)$$

to the numbers (8):

$$\begin{cases} \theta \cdot 2^r = \mathrm{K}^-_{1,\theta,r} - 1, \\ \mathrm{K}^-_{1,\theta,r+2} = 2\mathrm{K}^-_{1,\theta,r} - 3, \end{cases} (a) \quad \begin{cases} \theta \cdot 2^s = \mathrm{K}^-_{1,\theta,s} + 1, \\ \mathrm{K}^-_{1,\theta,s+1} = 2\mathrm{K}^-_{1,\theta,s} + 3, \end{cases} (b) \quad \text{and} \quad \begin{cases} \theta \cdot 2^s = \mathrm{K}^+_{1,\theta,s} - 1, \\ \mathrm{K}^+_{1,\theta,s+2} = 2\mathrm{K}^+_{1,\theta,s} - 3, \end{cases} (c) \quad \begin{cases} \theta \cdot 2^r = \mathrm{K}^+_{1,\theta,r} + 1, \\ \mathrm{K}^+_{1,\theta,r+2} = 2\mathrm{K}^+_{1,\theta,r} + 3. \end{cases} (d) \Rightarrow \begin{cases} \mathrm{K}^+_{1,\theta,r} = \mathrm{K}^-_{1,\theta,s}, \\ \mathrm{K}^-_{1,\theta,r} = \mathrm{K}^+_{1,\theta,s}. \end{cases} (f). \qquad (8)$$

We observe that transformations (6) - (8) lead to two groups of similar numbers (8f):

$$\mathrm{K}^{-(+)}_{a,\theta,s(r),} = \frac{1}{a}\left[\theta \cdot 2^{r(s)} - 1\right] \text{ (a) and } \mathrm{K}^{+(-)}_{a,\theta,r(s),} = \frac{1}{a}\left[\theta \cdot 2^{r(s)} + 1\right] \text{(b)}. \qquad (9)$$

Let's introduce the notations for the numbers (9): $\mathrm{K}^{-(+)}_{a,\theta,s(r),} = m_{a,\theta,r(s)}$ and $\mathrm{K}^{+(-)}_{a,\theta,r(s),} = p_{a,\theta,s(r)}$, and show that the numbers $m(p)_{a,\theta,r(s)}$ define nodes on the sequences $\theta \cdot 2^{r(s)}$, through which so-called the Jacobsthal tree ($JT$) branches in the $n \to \infty$ direction, forming $CT_{aq\pm 1}$ in the reverse $n \to 0$ direction. For this, we consider specific models of numbers (5) with parameter $a$.

**Case** $a = 1$. In this case, the numbers (5) are equals:

$$\mathrm{K}^{\pm}_{1,\theta,n} = \theta \cdot 2^n \pm (-1)^n \quad . \qquad (1.1)$$

We calculate the numbers (1.1) for $\theta = 1 \div 9$:

Table 1.1. Numbers $\mathrm{K}^{\pm}_{1,\theta,n}$

| | | | | | | | | | | |
|---|---|---|---|---|---|---|---|---|---|---|
| 9 | 10 | 17 | 37 | 71 | 145 | 287 | 577 | 1151 | $\mathrm{K}^+_{1,9,n}$ | Unknown (Unkn) |
| 7 | 8 | 13 | 29 | 55 | 113 | 223 | 449 | 895 | $\mathrm{K}^+_{1,7,n}$ | A321483 |
| 5 | 6/3=2 | 9/3=3 | 21/3=7 | 39/3=13 | 81=3/27 | 159/3=53 | 321/3=107 | 639/3=213 | $\mathrm{K}^+_{1,5,n}$ | Unkn/ A049887 |
| 3 | 4 | 5 | 13 | 23 | 49 | 95 | 193 | 769 | $\mathrm{K}^+_{1,3,n}$ | A140683 |
| 1 | 2 | 1 | 5 | 7 | 17 | 31 | 65 | 127 | $\mathrm{K}^+_{1,1,n}$ | **A014551** |
| $\theta_i$ \ $2^n$ | $2^0$ | $2^1$ | $2^2$ | $2^3$ | $2^4$ | $2^5$ | $2^6$ | $2^7$ | | OEIS [26] |
| 1 | 0/3=0 | 3/3=1 | 3/3=1 | 9/3=3 | 15/3=5 | 33/3=11 | 63/3=21 | 129/3=43 | $\mathrm{K}^-_{1,1,n}$ | A062510 / A001045 |
| 3 | 2 | 7 | 11 | 25 | 47 | 97 | 191 | **385** | $\mathrm{K}^-_{1,3,n}$ | A201630 |
| 5 | 4 | 11 | 19 | 41 | 79 | 161 | 319 | 641 | $\mathrm{K}^-_{1.5,n}$ | A321643 |
| 7 | 6/3=2 | 15/3=5 | 27/3=9 | 57/3=19 | 111/3=27 | 225/3=75 | 447/3=149 | 897/3=299 | $\mathrm{K}^-_{1,7,n,}$ | Unkn / Unkn |
| 9 | 8 | 19 | 35 | 73 | 143 | 289 | 575 | 1153 | $\mathrm{K}^-_{1,9,n}$ | Unkn |

For this numbers satisfy the recurrence relation

$$\mathrm{K}^{\pm}_{1,\theta,n+2} = \mathrm{K}^{\pm}_{1,\theta,n+1} + 2\mathrm{K}^{\pm}_{1,\theta,n}, \qquad (1.2)$$

with initial conditions

$$\begin{cases} \mathrm{K}^+_{1,\theta,0} = \theta + 1, \\ \mathrm{K}^+_{1,\theta,0} + \theta + 2 = \mathrm{K}^+_{1,\theta,1} \end{cases} \text{and} \quad \begin{cases} \mathrm{K}^-_{1,\theta,0} = \theta + 1, \\ \mathrm{K}^-_{1,\theta,0} + \theta - 2 = \mathrm{K}^-_{1,\theta,1} \end{cases}. \qquad (1.3)$$

For two adjacent numbers $\mathrm{K}^{\pm}_{1,\theta,n}$ in a row, the formula holds true:

$$\mathrm{K}^{\pm}_{1,\theta,n+1} = 2\mathrm{K}^{\pm}_{1,\theta,n} \mp 3. \qquad (1.4)$$

Rules (1.4) are analogous to rules of type $4x \pm 1$ in conjectures $3q \pm 1$.

The multiples of three $K^{-}_{1,1,n}$ numbers are known as recurrent Luke-Jacobsthal ($K^{+}_{1,1,n}$) numbers, and multiples of three numbers $K^{-}_{1,1,n}$ are known as recurrent Jacobsthal numbers [18,23 -25]. Numbers with index $\theta = 1$ of the type $2^n - 1$ are also known as Mersenne numbers for positive integer $n$, and of type $2^m + 1$ as Fermat numbers for non-negative integer $m = 2^n$. Therefore, with the index $\theta$, the numbers $K^{\pm}_{1,1,n}$ can still be structured as follows [9,10]. Numbers $K^{+}_{1,\theta,n}$, for which the condition for $\theta$ is fulfilled:

$$\frac{\theta - 1}{3} = Integer(I) \quad for \quad K^{-}_{1,\theta,n} \quad and \quad \frac{\theta + 1}{3} = I \quad for \quad K^{+}_{1,\theta,n}, \qquad (1.5)$$

multiples of three For example, if $\theta = 113$, then $\frac{113+1}{3} = 38$. Therefore, the multiples of three numbers are the numbers $K^{+}_{1,113,n}$: 114, 225, 453, 903, 1809,....

On the basis of (9), for the numbers $m_{1,\theta,r(s)}$ and $p_{1,\theta,r(s)}$ we write the equalities

$$\theta \cdot 2^{s(r)} = m_{1,\theta,s(r)} - 1 \quad \text{(a)} \qquad \text{and} \qquad \theta \cdot 2^{r(s)} = p_{1,\theta,r(s)} + 1 \quad \text{(b)}, \qquad (1.6)$$

and calculate them for $\theta = 1 \div 11$:

| | | | | | | | | | | | |
|---|---|---|---|---|---|---|---|---|---|---|---|
| 11 | 12 | 23 | 45 | 89 | 177 | 353 | 705 | 1409 | 2817 | A083683 | $p_{1,11,n}$ |
| 9 | 10 | 19 | 37 | 73 | 145 | 289 | 577 | 1153 | 2305 | A083705 | $p_{1,9,n}$ |
| 7 | 8 | 15 | 29 | 57 | 113 | 225 | 449 | 897 | 1793 | A083686 | $p_{1,7,n}$ |
| 5 | 6 | 11 | 21 | 41 | 81 | 161 | 321 | 641 | 1281 | A000383 | $p_{1,5,n}$ |
| 3 | 4 | 7 | 13 | 25 | 49 | 97 | 193 | 385 | 769 | A181565 | $p_{1,3,n}$ |
| 1 | 2 | 3 | 5 | 9 | 17 | 33 | 65 | 129 | 257 | **A000051** | $p_{1,1,n}$ |
| $\theta_i \backslash 2^n$ | $2^0$ | $2^1$ | $2^2$ | $2^3$ | $2^4$ | $2^5$ | $2^6$ | $2^7$ | $2^7$ | *OEIS* [26] | |
| 1 | 0 | 1 | 3 | 7 | 15 | 31 | 63 | 127 | 255 | A000225 | $m_{1,1,n}$ |
| 3 | 2 | 5 | 11 | 23 | 47 | 95 | 191 | 383 | 767 | A083329 | $m_{1,3,n}$ |
| 5 | 4 | 9 | 19 | 39 | 79 | 159 | 319 | 639 | 1279 | A153894 | $m_{1,5,n}$ |
| 7 | 6 | 13 | 27 | 55 | 111 | 223 | 447 | 895 | 1791 | A086224 | $m_{1,7,n}$ |
| 9 | 8 | 17 | 35 | 71 | 143 | 287 | 575 | 1151 | 2303 | A052996 | $m_{1,9,n}$ |
| 11 | 10 | 21 | 43 | 87 | 175 | 703 | 351 | 1407 | 2815 | A086225 | $m_{1,11,n}$ |

(1.7)

For adjacent numbers $m(p)_{1,\theta,n}$ in the rows, the recurrence formulas are true:

$$m_{1,\theta,n+1} = 2m_{1,\theta,n} + 1 \quad \text{and} \quad p_{1,\theta,n+1} = 2p_{1,\theta,n} - 1, \qquad (1.8)$$

with the initial conditions $m_{1,\theta,0}, p_{1,\theta,0}$, and for the numbers in each column, the formulas are executed:

$$m_{1,i+1,n} = m_{1,i,n} + 2^{n+1} \quad \text{and} \quad p_{1,i+1,n} = p_{1,i,n} + 2^{n+1}. \qquad (1.9)$$

The interval between recurrent numbers (nodes on the sequences $\theta \cdot 2^n$) $m(p)_{1,\theta,n}$ (1.8) is equal to:

$$m(p)_{1,\theta,n+1} - m(p)_{1,\theta,n} = m(p)_{1,\theta,n} \pm 1, \qquad (1.10)$$

which agrees with transformations (1.6).

As shown in Figure 1.1, according to the rules (1.6), the sequences $\theta \cdot 2^n$ branch out (indicated by arrows $\Rightarrow$) in the $n \to \infty$ direction from the nodes $m(p)_{1,\theta,r(s)}$ and form the *JT*, where the numbers of the nodes are the first elements $q_{odd} = \theta \cdot 2^0$ (or $\theta = m(p)_{a,\theta,0}$). Arrows $\to$ on the *JT* represent the $CT_{1q+1}$ for the numbers $q = 114$.

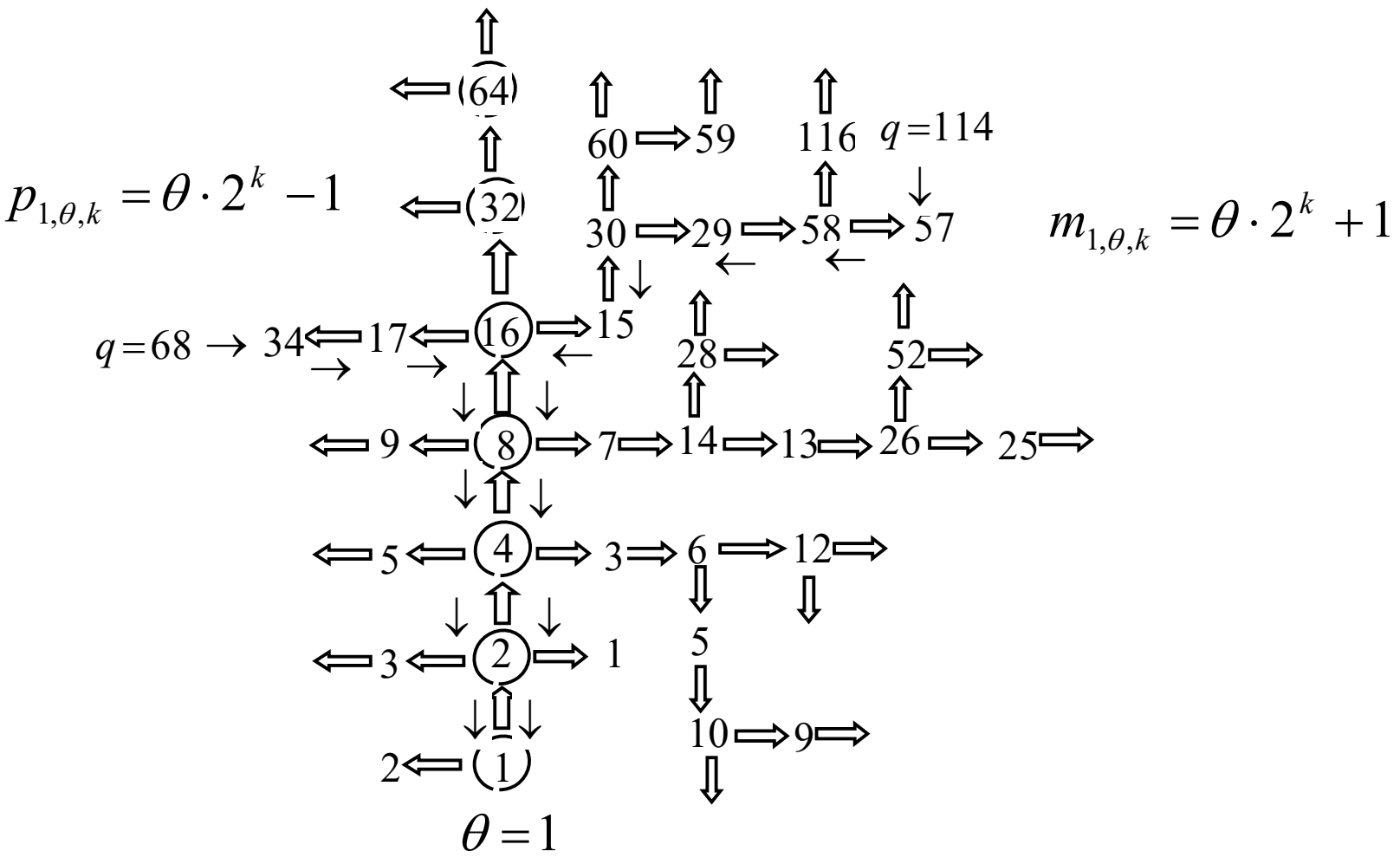

Figure 1.1. *Graphs of (→) $CT_{1q\pm1}$ conjecture and of (⇒) the JT*

In Figure 1.2, the $CT_{1q+1}$ of numbers $q = 34, 57$ are plotted in semi-logarithmic coordinates $N, \log(X(N))$, where v is the intermediate value of the transformed number, and $N$ is the iteration step. It can be seen that, unlike the conjecture $C_{q-1}$, the conjecture $C_{q+1}$ ends with a periodic cycle. Moreover, trajectories can be regarded as a superposition of tracks, each consisting of one conjectures $C_{q\pm1}$ and $k_{q/2}$. Then, between two odd numbers, the trajectories include an integer number of tracks. We demonstrate that for the model $C_{q\pm1}$, formula (1.11) is also valid:

$$\overbrace{\frac{\overbrace{\frac{\overbrace{\frac{\overbrace{q_{odd,end}\cdot 2^k - 1}^{track\ 1}}{a} 2^m - 1}^{track\ 2}}{a} 2^s - 1}^{track\ 3}}{a} \ldots 2^t - 1}^{track\ M} = q_{odd,initial} \qquad (1.11)$$

In formula (1.11), between two odd numbers $q_{odd,initial}$ and $q_{odd,end}$, , there are M tracks, where one track consists of one iteration $C_{q+1}$ and $k_{q/2}$ iterations $q/2$ (Figure 1.2). Therefore, in the $n \to \infty$ direction, in formula (1.11), the odd number $q_{odd,initial}$ is the initial value, and the odd number $q_{odd,end}$ is the final value, between which there are an integer number of tracks. For the $CT_{1q+1}$ of the number $q=57$ in Figure 1.2, $q_{odd,end}=1$, hence $k_{q/2}=4$, $m_{q/2}=1, s_{q/2}=1$. Then, the number $q_{odd,initial}$ equals $q_{odd,initial}=$ $=(((2^4-1)2-1)2-1)2=57$, which coincides with the initial number $q_{odd,initial}=57$.

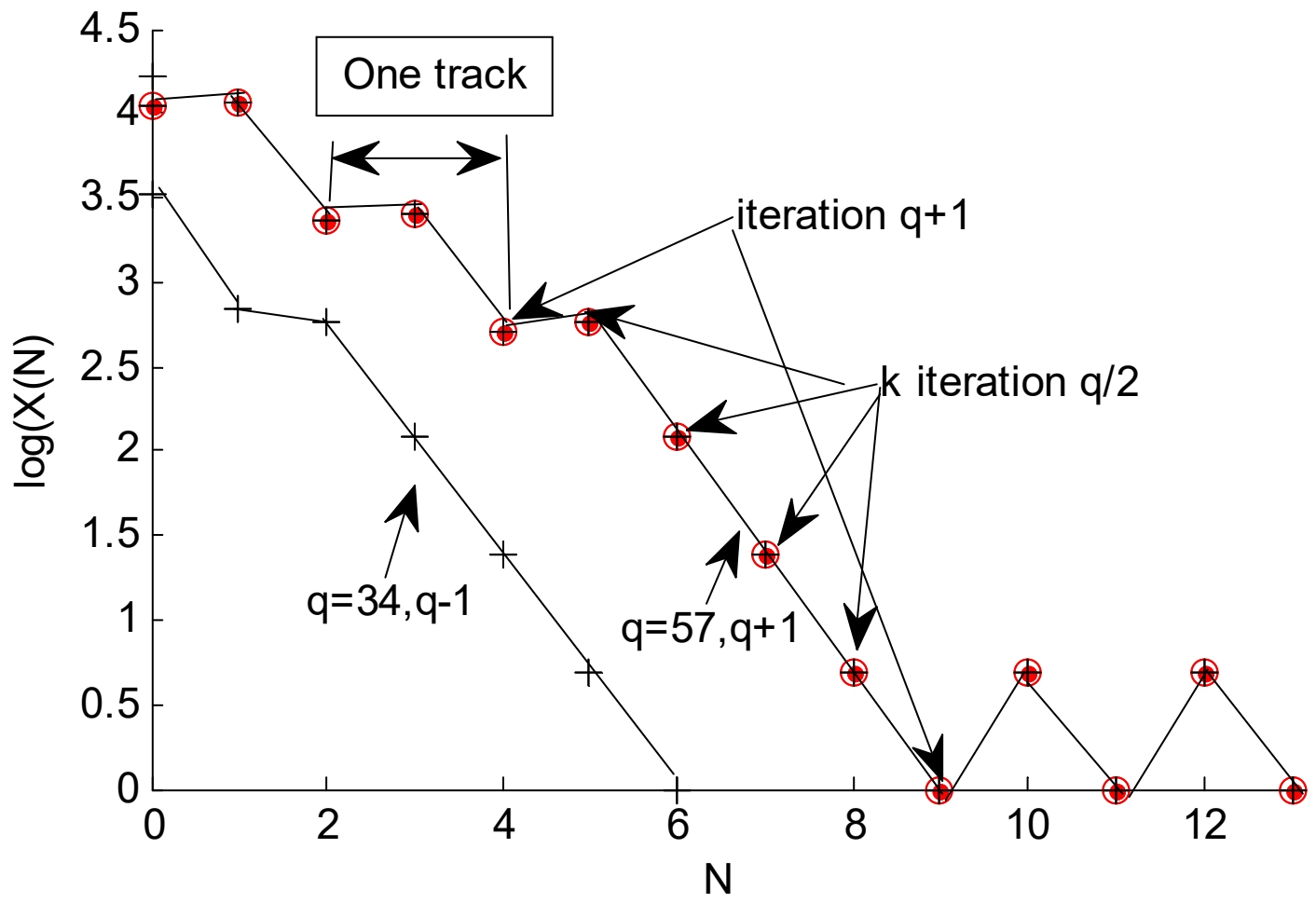

Figure 1.2

**Conclusion:** The regularities of transformations of recurrent numbers $K^{\pm}_{1,\theta,n}$ form the numbers $m(p)_{1,\theta,n}$ nodes on the sequences $\theta\cdot 2^n$ with a period $T_{a=1}=2^0$ and the rules of their transformations, from which the $JT$ branches out, forming $CT_{1q\pm1}$ in the reverse direction $n\to 0$. In each track

$$one\ \ track: \ \ q_{odd}\to q_{odd}\pm 1\to (q_{odd}\pm 1)/2^k = q^{\pm}_{odd,new}$$

| | | | | | | | | | | | | | | |
|---|---|---|---|---|---|---|---|---|---|---|---|---|---|---|
| $q^{-}_{odd,new}$ | 1 | 1 | 3 | 1 | 5 | 3 | 7 | 1 | 1 | 5 | 11 | 3 | 13 | … |
| | ↑ | ↑ | ↑ | ↑ | ↑ | ↑ | ↑ | ↑ | ↑ | ↑ | ↑ | ↑ | ↑ | … |
| $q_{odd}$ | 3 | 5 | 7 | 9 | 11 | 13 | 15 | 17 | 19 | 21 | 23 | 25 | 27 | … |
| | ↓ | ↓ | ↓ | ↓ | ↓ | ↓ | ↓ | ↓ | ↓ | ↓ | ↓ | ↓ | ↓ | … |
| $q^{+}_{odd,new}$ | 1 | 3 | 5 | 7 | 9 | 11 | 13 | 15 | 17 | 19 | 21 | 23 | 25 | … |

of the $CT_{1q+1}$, the inqualit $q_{odd}\ \langle\ q^{\pm}{}_{odd,new}$ is fulfilled, so for all numbers $q\in\mathbb{N}$ a unit element is reached.

According (3), for the numbers $q_{even}\neq\theta\cdot 2^n$, the number $k_{q/2}\leq 2$. Therefore, if in the track of the periodic cycle $k_{q/2}=1$, then the equality $q_{odd}=\dfrac{q_{odd}\pm 1}{2^1}$ holds only for $q_{odd}=\pm 1$. Thus, the periodic cycle with the minimum odd number $q_{odd}=+1$ can arise only in the $CT_{1q+1}$ conjecture. Periodic cycles with two $k_{q/2}=2$ iterations in the $CT_{1q+1}$ conjecture do not occur.

**Case 2.** $a = 3$. Selection of integers by formula

$$K^{\pm}_{3,\theta,n} = \frac{1}{3}\left[\theta \cdot 2^n \pm (-1)^n\right]. \qquad (2.1)$$

leads to the structuring of recurrent numbers in the form

$$K^{+}_{3,\theta,s(r)} = fractal, \quad K^{-}_{3,\theta,r(s)} = I \qquad \text{or} \qquad K^{+}_{3,\theta,s(r)} = I, \quad K^{-}_{3,\theta,r(s)} = fractal. \qquad (2.2)$$

and calculate them for $\theta = 1 \div 13$:

| | | | | | | | | | | | | |
|---|---|---|---|---|---|---|---|---|---|---|---|---|
| 3 | 13 | 4 | 9 | 17 | 35 | 69 | 139 | 277 | 555 | 1109 | $K^{-}_{3,13,n}$ | *Unkn* |
| | 11 | - | - | - | - | - | - | - | - | - | $K^{-}_{3,11,n}$ | |
| 2 | 7 | 2 | 5 | 9 | 19 | 37 | 75 | 149 | 299 | 597 | $K^{-}_{3,7,n}$ | *A062092* |
| | 5 | - | - | - | - | - | - | - | - | - | $K^{-}_{3,5,n}$ | |
| 1 | 1 | 0 | 1 | 1 | 3 | 5 | 11 | 21 | 43 | **85** | $K^{-}_{3,1,n}$ | *A001045* |
| *i* | $\theta_i \backslash 2^n$ | $2^0$ | $2^1$ | $2^2$ | $2^3$ | $2^4$ | $2^5$ | $2^6$ | $2^7$ | $2^8$ | | *OEIS[26]* |
| | 1 | - | - | - | - | - | - | - | - | - | $K^{+}_{3,1,n}$ | |
| 1 | 5 | 2 | 3 | 7 | 13 | 27 | 53 | 107 | 213 | 427 | $K^{+}_{3,5,n}$ | *A048573* |
| | 7 | - | - | - | - | - | - | - | - | - | $K^{+}_{3,7,n}$ | |
| 2 | 11 | 4 | 7 | 15 | 29 | 59 | 117 | 235 | 469 | 939 | $K^{+}_{3,11,n,}$ | *Unkn* |
| | 13 | - | - | - | - | - | - | - | - | - | $K^{+}_{3,13,n}$ | |

(2.3)

The numbers $K^{\pm}_{3,\theta,n}$ are the recurrence formula holds

$$K^{\pm}_{3,\theta,n+2} = K^{\pm}_{3,\theta,n+1} + 2K^{\pm}_{3,\theta,n}, \qquad (2.4)$$

with initial conditions :

$$If \quad \frac{\theta \pm 1}{3} = I \quad then \quad K^{\pm}_{3,\theta,0} = \frac{\theta \pm 1}{3} \quad and \quad K^{\pm}_{3,\theta,1} = \theta - K^{\pm}_{3,\theta,0} \quad . \qquad (2.5)$$

In rows (2.3), for two adjacent numbers $K^{\pm}_{3,\theta,n}$, the formula holds true:

$$K^{\pm}_{3,\theta,n+1} = 2K^{\pm}_{3,\theta,n} \mp 1, \qquad (2.6)$$

and in the column

$$K^{\pm}_{3,i+1,n} = K^{\pm}_{3,i,n} + 2^{n+1}, \quad \text{where } i \in \mathbb{N}. \qquad (2.7)$$

Taking (9) into account, we write down the formulas for the numbers $m(p)_{3,\theta,n}$:

$$m_{3,\theta,n} = \frac{1}{3}\left[\theta \cdot 2^n - 1\right] \quad (a) \qquad \text{and} \qquad p_{3,\theta,n} = \frac{1}{3}\left[\theta \cdot 2^n + 1\right] \quad (b), \qquad (2.8)$$

which are calculated for $\theta = 1 \div 13$:

| | | | | | | | | | | | |
|---|---|---|---|---|---|---|---|---|---|---|---|
| 5 | 13 | | 9 | | 34 | | 139 | | 555 | | $p_{3,13,n}$ |
| 4 | 11 | 4 | | 15 | | 59 | | 235 | | 939 | $p_{3,11,n}$ |
| 3 | 7 | | 5 | | 19 | | 75 | | 299 | | $p_{3,7,n}$ |
| 2 | 5 | 2 | | 7 | | 27 | | 107 | | 427 | $p_{3,5,n}$ |
| 1 | 1 | | 1, k=1 | | 11, k=2 | | 43, k=3 | | 171, k=4 | | $p_{3,1,n}$ |
| *i* | $\theta_i \backslash 2^n$ | $2^0$ | $2^1$ | $2^2$ | $2^3$ | $2^4$ | $2^5$ | $2^6$ | $2^7$ | $2^8$ | |
| 1 | 1 | 0, k=1 | | 1, k=2 | | 5, k=3 | | 21, k=4 | | 85, k=5 | $m_{3,1,n}$ |
| 2 | 5 | | 3 | | 13 | | 53 | | 213 | | $m_{3,5,n}$ |
| 3 | 7 | 2 | | 9 | | 37 | | 149 | | 597 | $m_{1,7,n}$ |
| 4 | 11 | | 7 | | 29 | | 117 | | 469 | | $m_{3,11,n}$ |
| 5 | 13 | 4 | | 17 | | 69 | | 277 | | 1109 | $m_{3,13,n}$ |

(2.9)

For the numbers $m(p)_{3,\theta,n}$ in the rows, the formulas are true:

$$m_{3,\theta,k+1} = 4m_{3,\theta,k} + 1 \quad \text{and} \quad p_{3,\theta,k+1} = 4p_{3,\theta,k} - 1. \qquad (2.10)$$

and in the columns:

$$m_{3,i+2,n} = m_{3,i,n} + 2^{n+1} \quad \text{and} \quad p_{3,i+2,n} = p_{3,i,n} + 2^{n+1}. \qquad (2.11)$$

According (2.8), we have:

$$m(p)_{3,\theta,k+1} - m(p)_{3,\theta,k} = 4m(p)_{3,\theta,k} \pm 1 - m(p)_{3,\theta,k} = 3m(p)_{3,\theta,k} \pm 1. \qquad (2.12)$$

Formula (2.12) agrees with formulas (2.8).

Integers

$$m_{3,\theta_{1,5},r(s)} = \frac{\theta_{1,5} \cdot 2^{r(s)} - 1}{3} = I \text{ (a)} \quad \text{and} \quad p_{3,\theta_{1,5},r(s)} = \frac{\theta_{1,5} \cdot 2^{r(s)} + 1}{3} = I \text{ (b)} \qquad (2.13)$$

correspond to nodes (see also *Appendix 3*) of sequences $\theta \cdot 2^n$ from which sequences with different indices $\theta$ grow in the direction $n \to \infty$. This is how the *JT* is formed. If

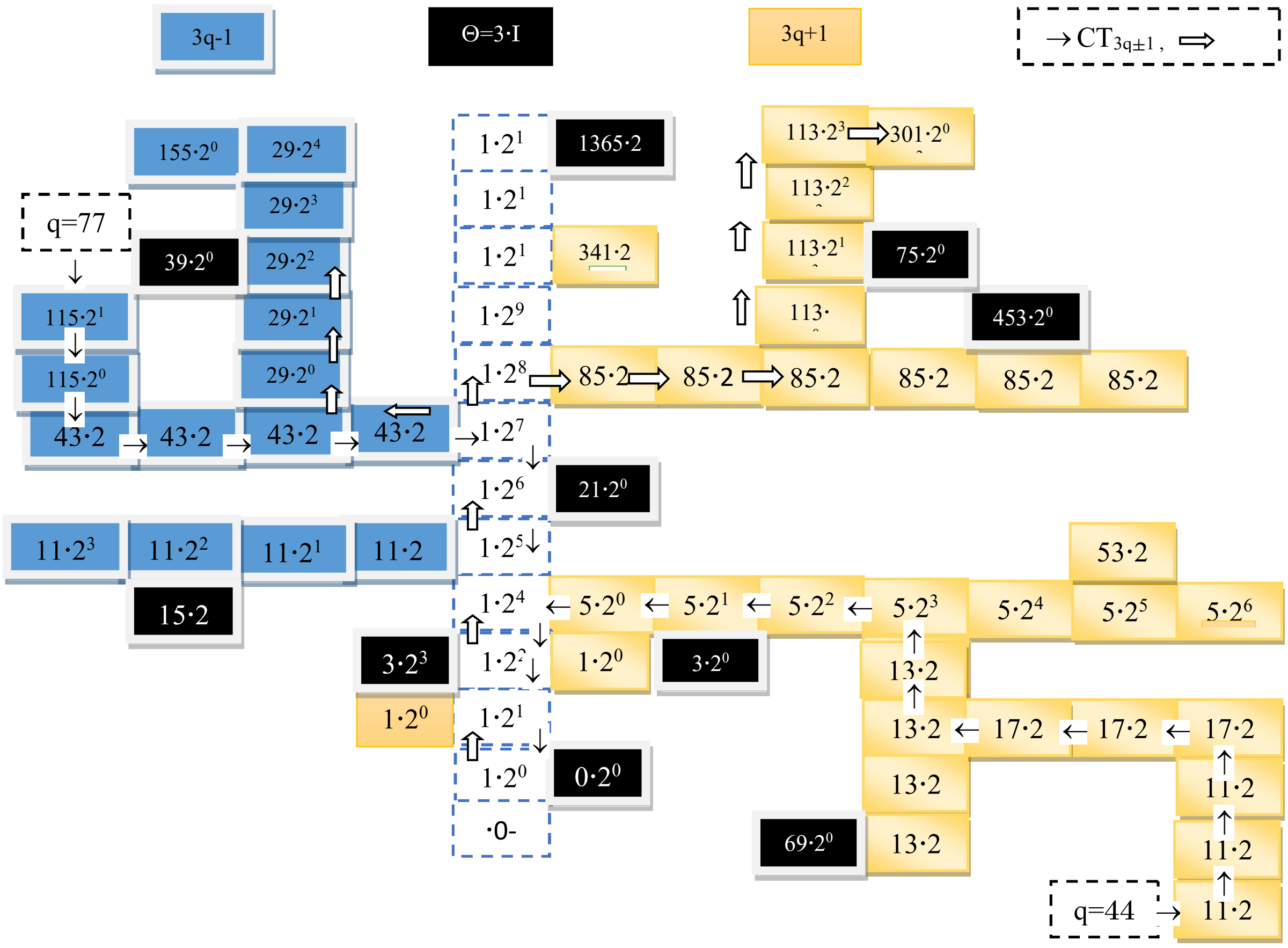

**Figure 2.1.** *Graphs of* ($\rightarrow$) $CT_{3q\pm1}$ *conjecture and of* ($\Rightarrow$) *the JT*

the index of the sequence $\theta$ is a multiple of three

$$\theta = \theta_3 = I \cdot 3, \qquad (2.14)$$

then of the numbers $m(p)_{3,\theta_{1,5},n}$ of the fractions and the *JT* do not formed (see also *Appendix 4*).

An illustration of cycles $cycle_{1\leftrightarrow1}^{3q\pm1}$ the completion of their trajectories $CT_{3q\pm1}$ in the direction of the formation of single elements is shown in Figure 2a,b. For the conjecture $3q-1$ isolated from the numbers $1\cdot2^n$ periodic cycles with minimum values of oscillations $q_{odd,\min}=5$ and $q_{odd,\min}=17$ are observed (Figure 2.2a-c,f). The number of branches on the circle is determined by the fulfillment of condition (2.13). As shown in Figure 2.3, isolated from the numbers $1\cdot2^n$ are the cycles $cycle_{5\leftrightarrow5}^{3q-1}$ and $cycle_{17\leftrightarrow17}^{3q-1}$, caused by periodic oscillations in the *JT*. It is shown also in *Appendix 5,* the set of numbers $q\in\mathbb{N}$ is divided in approximately the same proportions: $cycle_{1\leftrightarrow1}^{3q-1}$ (33%), $cycle_{5\leftrightarrow5}^{3q-1}$ (32%) and $cycle_{17\leftrightarrow17}^{3q-1}$ (35%) from the first 300,000 numbers taken.

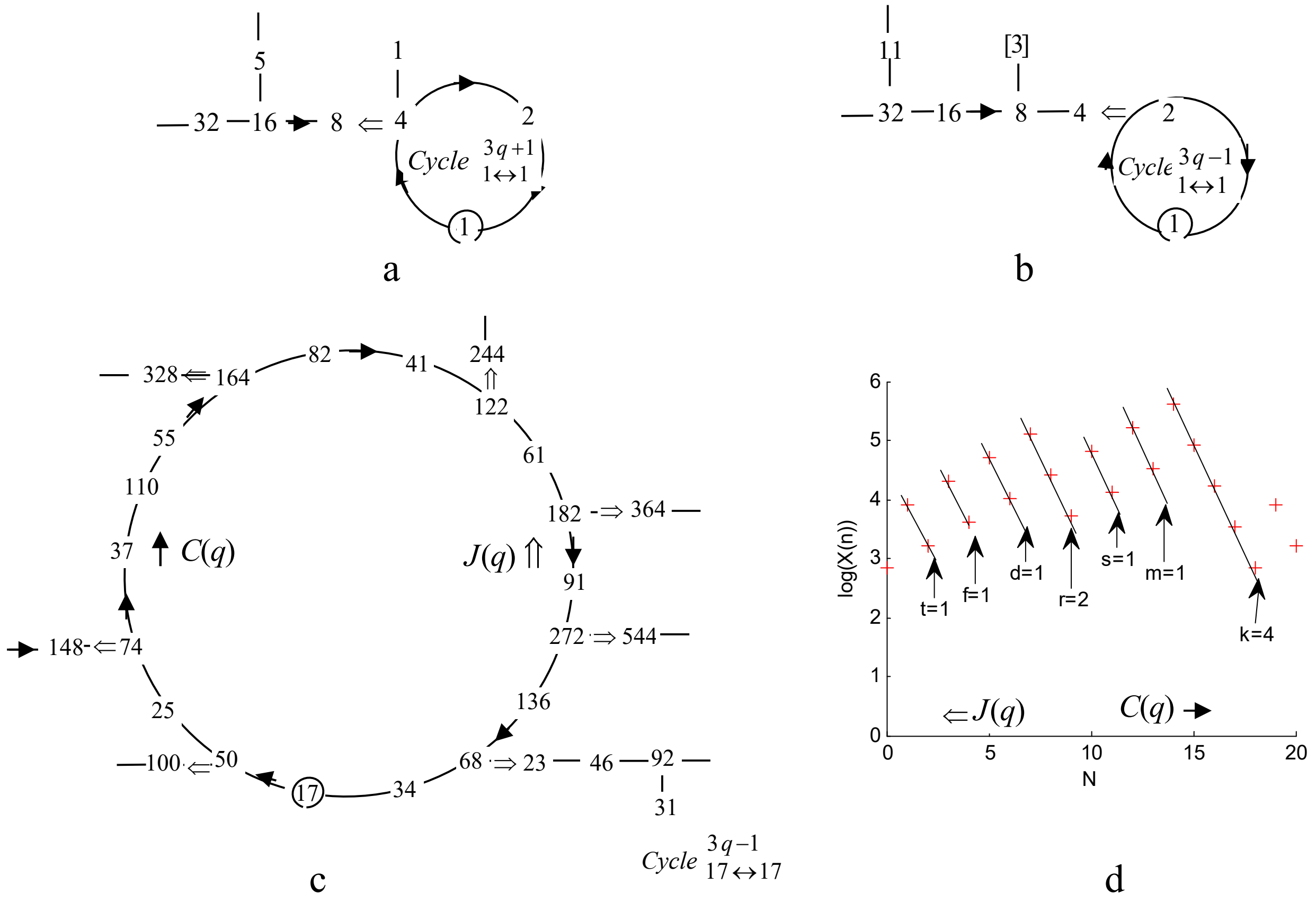

**Figure 2.2.** *The Jacobsthal and the Collatz graphs (a,b,c) and the conjecture q=17 in the semi-logarithmic scale* $N,\log q$ *(d)*

In Figure 2.1, arrows of type $\Rightarrow$ show a fragment of the *JT*. In the $n\to\infty$ direction to the right relative to the sequence $1\cdot2^n$, the tree branches according to the rule (2.8a), and in the direction to the left relative to the numbers $1\cdot2^n$ , the tree branches according to the rule (2.8b). Here arrows of type $\rightarrow$ show a two $CT_{3q\pm1}$ of numbers $q=44\,(3q+1)$ and $q=77\,(3q-1)$ (see also *Appendix 2*).

Let's investigate the conditions for the formation of periodic cycles. A loop occurs when an intermediate value of the converted number is repeated. If the period of the cycle is formed from one conjecture $3q\pm1$ and $k_{q/2}$ iterations $q/2$ (one track), then the equality is holds:

$$\frac{\overbrace{q_{odd,k}\cdot 2^k \mp 1}^{track}}{3}=q_{odd,k}=I,\ k=k_{q/2} \qquad (2.15)$$

Then for the conjecture $3q-1$ from (2.15) we have, that $q_{odd,1}=1$ at $k_{q/2}$=1, that is, the $cycle_{1\leftrightarrow1}^{3q-1}$ is formed by one iteration of $q/2$ . For conjecture $3q+1$, the value of $q_{odd,2}=1$ at $k_{q/2}$=2, i.e. the loop is formed from two iterations of $q/2$. At $k_{q/2}$>2, $q_{k,\min}\neq I$.

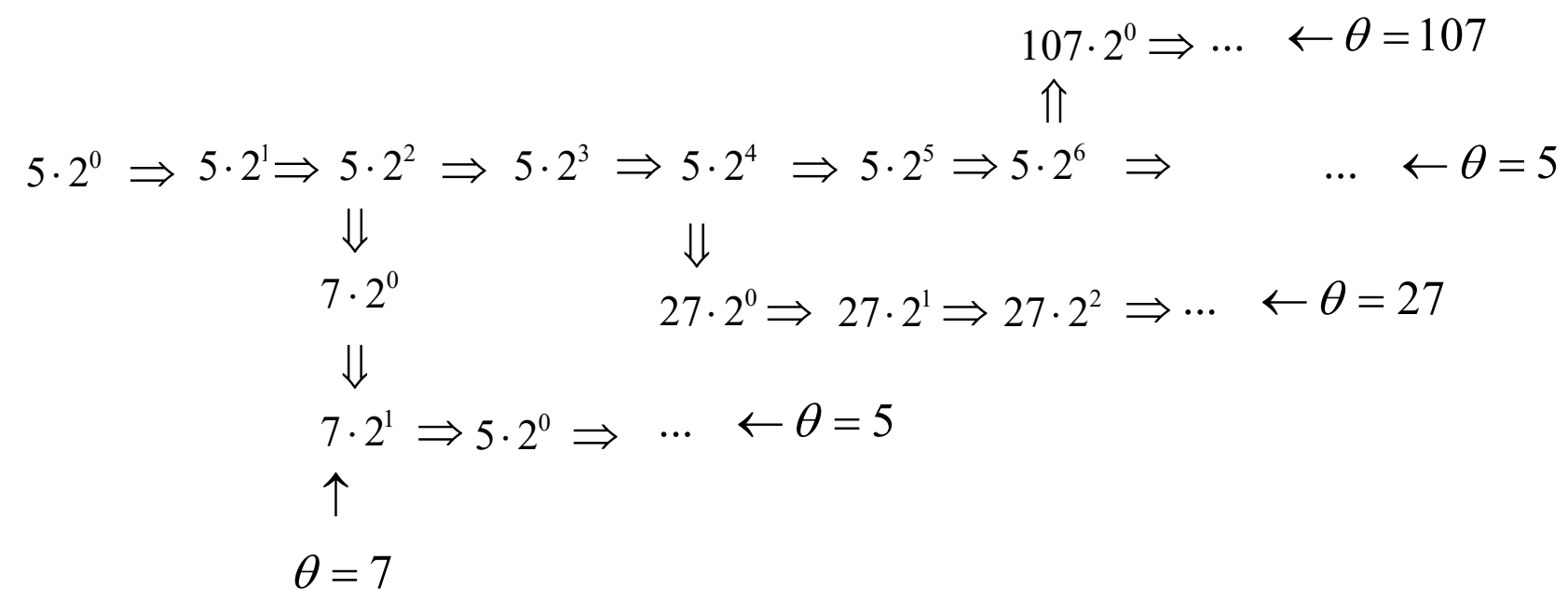

**Figure 2.3** *The graph of isolated cycle* of the *JT for $C_{3q-1}$ conjecture*

The cycle $cycle_{5\leftrightarrow5}^{3q-1}$ consists of two tracks (2.15). In this case, the equality holds:

$$\frac{\overbrace{\frac{\overbrace{q_{odd,m,k}\cdot 2^k \mp 1}^{track}}{3}2^m \mp 1}^{track}}{3}=q_{odd,m,k}\Rightarrow q_{odd,k,m}=\frac{\pm(3+2^m)}{2^{k+m}-3^2},\quad k(m)=k_{q/2}(m_{q/2}) \qquad (2.16)$$

Therefore, according to (2.16), for the conjecture $3q+1$ we have one cycle with the parameter $q_{odd,2,2}=1$ at $k_{q/2}$=2 and $m_{q/2}$=2, while for the conjecture $3q-1$ we have three cycles with the parameters $q_{odd,1,1}=1$ at $k_{q/2}$=1, $m_{q/2}$ =1, and $q_{odd,2,1}=5$ at $k_{q/2}$=2, $m_{q/2}$=1, and $q_{odd,1,2}=7$ at $k_{q/2}$=1, $m_{q/2}$ =2 .

The $cycle_{17\leftrightarrow17}^{3q-1}$ is formed from seven steps (2.16) (Figure 2.2d). Therefore, the equality holds for it:

$$\overbrace{\frac{\overbrace{\frac{\overbrace{\frac{\overbrace{q_{odd,k,m,\ldots,t}\cdot 2^k-1}^{track1}}{3}2^m-1}^{track2}}{3}2^s-1}^{track3}}{3}\ldots 2^t-1}^{track7}=q_{odd,k,m,\ldots,t}\quad (a)\ \Rightarrow \qquad (2.17)$$

$$\Rightarrow q_{odd,k,m,\ldots,t}=\frac{\pm(3^6+3^5\cdot2^t+3^4\cdot2^{t+f}+3^3\cdot2^{d+f+t}+3^2\cdot2^{r+d+f+t}+3^1\cdot2^{s+r+d+f+t}+2^{m+s+r+t+d+f+t})}{-3^7+2^{k+m+s+r+t+t+d+f}},\quad (b)$$

where $k,...,d = k_{q/2},...,t_{q/2}$ . According to (2.17), the parameter $q_{odd,4,1,1,2,1,1,1} = \frac{2363}{139} = 17$ at $k_{q/2} = 4$ , $m_{q/2} = 1$ , $s_{q/2} = 1$ , $r_{q/2} = 2$ , $d_{q/2} = 1$ , $f_{q/2} = 1$ , $t_{q/2} = 1$ (Figure 2.2d). At other integer $k_{q/2},...,t_{q/2}$ , the values $q_{odd,k,...,t} \neq I$ or negative.

**Conclusion**: Conjecture $CT_{3q\pm1}$ is given the greatest attention, so we will give extended conclusions.

1. The regularities of transformations of recurrent numbers $K^{\pm}_{3,\theta,n}$ form the numbers $m(p)_{3,\theta,n}$ nodes on the sequences $\theta \cdot 2^n$ with a period $T_{a=3} = 2^1$ and the rules of their transformations, from which the *JT* branches out, forming of the $CT_{3q\pm1}$ in the reverse $n \to 0$ direction. In the first track $one\ \ track: q_{odd} \to 3q_{odd} - 1 \to \frac{3q_{odd} - 1}{2^v} = q^-_{odd,new}$ of the $CT_{3q\pm1}$, inquality of type $q_{odd} \langle q^-_{odd,new}$ and $q_{odd} \rangle q^-_{odd,new}$ are fulfilled. The distribution numbers $q^-_{odd,new}$

| $q^-_{odd,new}$ | 1 | 7 | 5 | 13 | 1 | 19 | 11 | 25 | 7 | 21 | 17 | 37 | 5 | … |
|---|---|---|---|---|---|---|---|---|---|---|---|---|---|---|
| | ↑ | ↑ | ↑ | ↑ | ↑ | ↑ | ↑ | ↑ | ↑ | ↑ | ↑ | ↑ | ↑ | … |
| $q_{odd}$ | 3 | 5 | 7 | 9 | 11 | 13 | 15 | 17 | 19 | 21 | 23 | 25 | 27 | … |

consistent with the regularities of the distribution of numbers in table A6-1.

2. Problems $C_{3q\pm1}$ of conjectures of natural numbers $q \in \mathbb{N}$ are formed regularities of transformations of parameterized $J_{\theta,n}$ recurring Jacobsthal numbers

$$m(p)_{3,\theta,n} = \left.\frac{\theta \cdot 2^n \pm (-1)^n}{3}\right|_{n \in N_{even}+N_{odd}} = I \Rightarrow \begin{cases} \frac{\theta \cdot 2^{N_{even}} \pm (-1)^{N_{even}}}{3} \Rightarrow \begin{cases} \frac{\theta \cdot 2^{N_{even}} + 1}{3} \Rightarrow C_{3q-1} \\ \frac{\theta \cdot 2^{N_{even}} - 1}{3} \Rightarrow C_{3q+1} \end{cases} \\ \frac{\theta \cdot 2^{N_{odd}} \pm (-1)^{N_{odd}}}{3} \Rightarrow \begin{cases} \frac{\theta \cdot 2^{N_{odd}} - 1}{3} \Rightarrow C_{3q+1} \\ \frac{\theta \cdot 2^{N_{odd}} + 1}{3} \Rightarrow C_{3q-1} \end{cases} \end{cases}$$

3. Doubling in the $n \to \infty$ direction of parameterized $\theta \neq 3 \cdot I = \theta_{1,5}$ Jacobsthal numbers $m(p)_{3,\theta,0} \cdot 2^k$ with further branching in the nodes $\frac{m(p)_{3,\theta,0} \cdot 2^k \mp 1}{3} = m(p)_{3,\theta_{new},0} \cdot 2^0$, forms a *JT*, on which in the reverse $n \to 0$ Collatz trajectories lead to the periodic cycle $cycle^{3q+1}_{1\leftrightarrow1}$ for the conjecture $C_{3q+1}$ and one of the cycles $cycle^{3q-1}_{1\leftrightarrow1}$, $cycle^{3q-1}_{5(7)\leftrightarrow5(7)}$ $and$ $cycle^{3q-1}_{17\leftrightarrow17}$ for the conjecture $C_{3q-1}$.

4. Parameterized $\theta = 3 \cdot I = \theta_3$ nodes do not branch.

5. Parameterized $\theta \neq 3 \cdot I = \theta_{1,5}$ numbers are structured in binary on groups: $\theta_1 = \frac{\theta - 1}{3} = I$ and $\theta_5 = \frac{\theta + 1}{3} = I$.

6. Nodes of sequences $\theta_{1,5}\cdot 2^{n}$ form periodic structures of type $\frac{\theta_5+1}{3}=Integer$

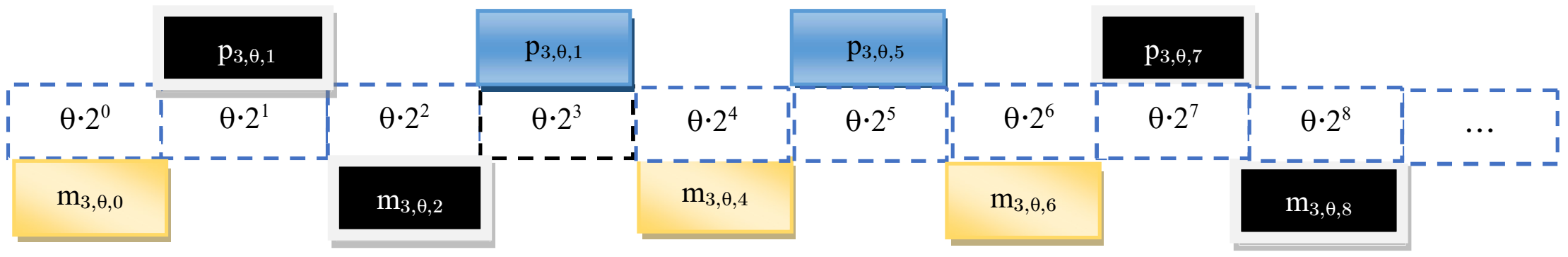

and type $\frac{\theta_1-1}{3}=Integer$

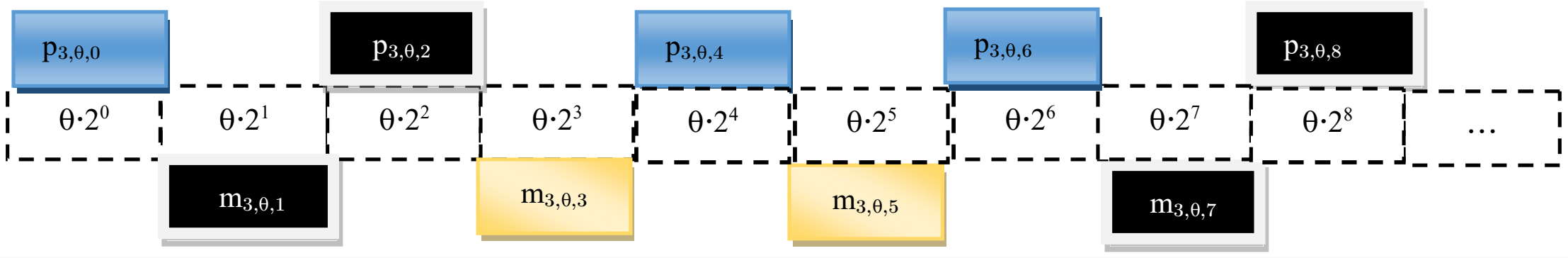

from which the *JT* branches. The period $T_n$ between two adjacent nodes, does not depend on the power $n$ and grows proportionally to $\theta:T_n=\theta\left(2^{n+2}-2^{n}\right)=\theta\cdot 2^{2}$.

7. By periodic cycles $cycle^{3q-1}_{1\leftrightarrow 1}$, $cycle^{3q-1}_{5(7)\leftrightarrow 5(7)}$ and $cycle^{3q-1}_{17\leftrightarrow 17}$, completion of the Collatz conjectures of numbers $q\in\mathbb{N}$, the set $q\in\mathbb{N}$ tends to a uniform distribution. The minimum odd numbers in cycles $cycle^{3q-1}_{1\leftrightarrow 1}$, $cycle^{3q-1}_{5(7)\leftrightarrow 5(7)}$ and $cycle^{3q-1}_{17\leftrightarrow 17}$ are localized in the interval of small values. Therefore, for the conjecture $C_{3q+1}$ the absence of similar cycles with values that are not equal to one is a "strong" argument in favor of the Collatz hypothesis.

8. In the interval of odd values ($q_{odd,initial},q_{odd,end}$), the trajectory $CT_{3q\pm 1}$ is formed by the integer $N$ number of tracks , so the equality holds:

$$\overbrace{\overbrace{\overbrace{\overbrace{\frac{\frac{\frac{q_{odd\,end}\cdot 2^{k}-1}{3}2^{m}-1}{3}2^{s}-1}{3}}^{track1}}^{track2}}^{track3}}^{track\,M}..2^{s}-1=q_{odd\,initial}\quad k,...=k_{q/2},...$$

9. According to (3), for the numbers $q_{even}\neq\theta\cdot 2^{n}$, the number $k_{q/2}\leq 2$. Therefore, if in the track of the periodic cycle $k_{q/2}=1$, then the equality $q_{odd}=\frac{3q_{odd}\pm 1}{2^{1}}$ holds only for $q_{odd}=\pm 1$. Thus, the periodic cycle with the minimum odd number $q_{odd}=+1$ can arise only in the $C_{3q-1}$ conjecture, as $1\rightarrow 2\rightarrow 1\rightarrow\ldots$ . If $k=2$ then $q_{odd}=\frac{q_{odd}\cdot 2^{n}-1}{4}\Rightarrow q_{odd}=$ -1. Therefore periodic cycles with two $k=2$ iterations in $C_{3q-1}$ conjecture do not occur.

If periodic cycle consists of two tracks with $k_{q/2}$ and $m_{q/2}$ numbers iterations, then $q_{odd}=\frac{\left[\frac{3q_{odd}\pm 1}{2^{k_{q/2}}}\right]3\pm 1}{2^{m_{q/2}}} \Rightarrow q_{odd}=\pm\frac{3+2^{k_{q/2}}}{2^{m_{q/2}+k_{q/2}}-3^2} \Rightarrow \begin{matrix} k_{q/2} & 1 & 2 & 1 & 2 \\ m_{q/2} & 1 & 2 & 2 & 1 \\ q_{odd} & \mp 1 & \pm 1 & \mp 5 & \mp 7 \end{matrix}$. Therefore, for the $C_{3q+1}$ conjecture , a periodic cycle with the minimum value of an odd number $q_{odd}=+1$ with parameters $k_{q/2}$=2 and $m_{q/2}$=2 is allowed, as two identical tracks of the type $1\to 4\to 2\to 1\to\ldots$. For the conjecture $C_{3q-1}$, periodic cycles with parameters $k_{q/2}$=1 and $m_{q/2}$=1, and with parameters $k_{q/2}$=2 and $m_{q/2}$=2 with the minimum value of an odd number $q_{odd}=+5$. are allowed. The known periodic cycle $cycle_{1\leftrightarrow 1}$ in the conjecture $C_{3q-1}$ is calculated similarly as a superposition of a larger number of tracks (Figure 2.2d)

**Case 3.** $a=5$. In this model, formula (5) is written in the form:

$$\mathrm{K}^{\pm}_{5,\theta,n}=\frac{1}{5}\left[\theta\cdot 2^n \pm(-1)^n\right]. \tag{3.1}$$

As follows from (3.2)

| θ$_i$ \2$^n$ | $2^0$ | $2^1$ | $2^2$ | $2^3$ | $2^4$ | $2^5$ | $2^6$ | $2^7$ | $2^8$ | | OEIS[26] |
|---|---|---|---|---|---|---|---|---|---|---|---|
| 1 | - | - | 1 | - | 3 | - | 13 | - | **51** | $\mathrm{K}^{\pm}_{3,1,n}$ | *A015521* |
| 3 | - | 1 | - | 5 | - | 19 | - | 77 | - | $\mathrm{K}^{\pm}_{3,5,n}$ | *A108981* |
| 7 | - | 3 | - | 11 | - | 45 | - | 179 | - | $\mathrm{K}^{\pm}_{3,7,n}$ | *Unkn* |
| 11 | 2 | - | 9 | - | 35 | - | 141 | - | 563 | $\mathrm{K}^{\pm}_{3,11,n,}$ | *Unkn* |
| 13 | - | 5 | - | 21 | - | 83 | - | 333 | - | $\mathrm{K}^{\pm}_{3,13,n}$ | *Unkn* |

(3.2)

the recurrence relation

$$\mathrm{K}^{\pm}_{5,\theta,n+2}=3\mathrm{K}^{\pm}_{5,\theta,n+1}+4\mathrm{K}^{\pm}_{5,\theta,n} \tag{3.3}$$

holds true for numbers $\mathrm{K}^{\pm}_{5,\theta,n}$. The analysis of the numbers bisected by

$$\mathrm{K}^{-}_{5,\theta,n}=\frac{1}{5}\left[\theta\cdot 2^n-(-1)^n\right]\text{(a)} \quad \text{and} \quad \mathrm{K}^{+}_{5,\theta,m}=\frac{1}{5}\left[\theta\cdot 2^n+(-1)^n\right]\text{(b)}, \tag{3.4}$$

signs (3.5)

| | | | | | | | | | | | | |
|---|---|---|---|---|---|---|---|---|---|---|---|---|
| 3 | 13 | - | - | - | 21 | - | - | - | 333 | - | $\mathrm{K}^{+}_{3,13,n}$ | *Unkn* |
| | 11 | 2 | - | - | - | 35 | - | - | - | 563 | $\mathrm{K}^{+}_{3,11,n}$ | *Unkn* |
| 2 | 7 | - | 3 | - | - | - | 45 | - | - | - | $\mathrm{K}^{+}_{3,7,n}$ | *A266698* |
| | 3 | - | - | - | 5 | - | - | - | 77 | - | $\mathrm{K}^{+}_{3,5,n}$ | *A318236* |
| 1 | 1 | | - | - | - | 3 | - | - | - | **51** | $\mathrm{K}^{+}_{3,1,n}$ | *A182512* |
| *i* | θ$_i$ \2$^n$ | $2^0$ | $2^1$ | $2^2$ | $2^3$ | $2^4$ | $2^5$ | $2^6$ | $2^7$ | $2^8$ | | OEIS [26] |
| | 1 | - | - | 1 | - | - | - | 13 | - | - | $\mathrm{K}^{-}_{3,1,n}$ | *A299960* |
| 1 | 3 | - | 1 | - | - | - | 19 | - | - | - | $\mathrm{K}^{-}_{3,5,n}$ | *A182460* |
| | 7 | - | - | - | 11 | - | - | - | 179 | - | $\mathrm{K}^{-}_{3,7,n}$ | *Unkn* |
| 2 | 11 | - | - | 9 | - | - | - | 141 | - | - | $\mathrm{K}^{-}_{3,11,n,}$ | *Unkn* |
| | 13 | - | 5 | - | - | - | 83 | - | - | - | $\mathrm{K}^{-}_{3,13,n}$ | *Unkn* |

(3.5)

shows that the recurring formulas

$$K^{\pm}_{5,\theta,m+4} = 16K^{\pm}_{5,\theta,m} \pm 3 \ . \qquad (3.6)$$

are fulfilled for the numbers in the rows.

Taking (9) into account, we write down the formulas for the numbers $m(p)_{5,\theta,n}$:

$$m_{5,\theta,n} = \frac{1}{5}\left[\theta \cdot 2^n - 1\right] \ \text{(a)} \quad \text{and} \quad p_{5,\theta,n} = \frac{1}{5}\left[\theta \cdot 2^n + 1\right] \ \text{(b)}, \qquad (3.7)$$

and calculate for $\theta = 1 \div 11$:

| 11 | | | 9 | | | | 141 | | | | $p_{3,13,n}$ |
|---|---|---|---|---|---|---|---|---|---|---|---|
| 9 | 2 | | | | 29 | | | | 461 | | $p_{3,11,n}$ |
| 7 | | 3 | | | | [45] | | | | 717 | $p_{3,7,n}$ |
| 3 | | | | 5 | | | | 77 | | | $p_{3,5,n}$ |
| 1 | | | 1, k=1 | | | | 13, k=2 | | | | $p_{3,1,n}$ |
| $\theta_i$ \ $2^n$ | $2^0$ | $2^1$ | $2^2$ | $2^3$ | $2^4$ | $2^5$ | $2^6$ | $2^7$ | $2^8$ | $2^0$ | |
| 1 | [0], k=1 | | | | 3, k=2 | | | | 51, k=3 | | $m_{3,1,n}$ |
| 3 | | 1 | | | | 19 | | | | 307 | $m_{3,5,n}$ |
| 7 | | | | 11 | | | | 179 | | | $m_{1,7,n}$ |
| 9 | | | 7 | | | | [115] | | | | $m_{3,11,n}$ |
| 11 | 2 | | | | [35] | | | | 563 | | $m_{3,13,n}$ |

(3.8)

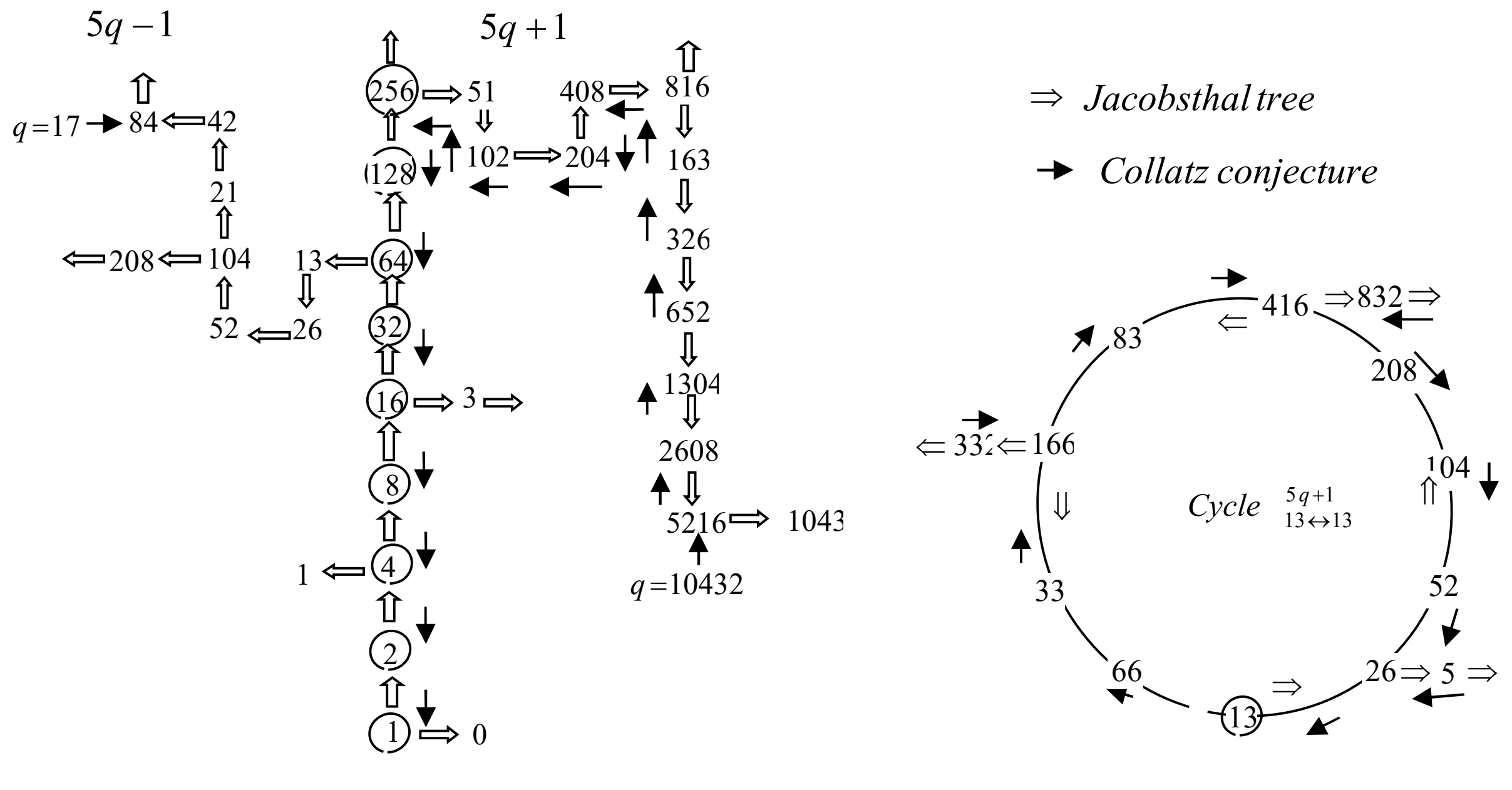

a b

Figure 3.1

In (3.8), square brackets highlight numbers that are multiples of five. For numbers $m(p)_{5,\theta,k}$ in lines, the conversion rules are followed

$$m_{5,\theta,k+1} = 16m_{5,\theta,k} + 3 \ \text{(a)} \quad \text{and} \quad p_{5,\theta,k+1} = 16p_{5,\theta,k} - 3 \ \text{(b)}, \qquad (3.9)$$

and for the numbers in the columns, formulas (1.9) have the form:

$$m_{5,\theta+1,n} = m_{5,\theta,n} + 2^{n+1} \quad \text{and} \quad p_{5,\theta+1,n} = p_{5,\theta,n} + 2^{n+1} . \qquad (3.10)$$

(3.9), we have:

$$m(p)_{5,\theta,k+1} - m(p)_{5,\theta,k} = 16m(p)_{5,\theta,k} \pm 3 - m(p)_{5,\theta,k} = 3(5m(p)_{5,\theta,k} \pm 1)\,. \qquad (3.11)$$

Formula (3.11) agrees with formulas (3.7).

An illustration of the application of the rule of numbering numbers $m(p)_{5,\theta,k}$ by index $k$ is shown for numbers $m(p)_{5,1,k}$. Rules (3.9) are analogous to rules of type $4x \pm 1$ in conjrctures $3q \pm 1$. Therefore, according to the rules (3.7), the numbers (3.1) form the nodes of the sequence $\theta \cdot 2^n$, from which other sequences $\vartheta \cdot 2^n$ branch out in the $n \to \infty$ direction of increasing power $n$, thus forming a tree of graphs, as shown in Figure 3.1a. An illustration of the trajectory of Collatz conjecture on the *JT* is given in Figure 3.1a for numbers $q = 17$ and $q = 10432$

**5x-1** **5x+1**

| | | | | | | | q | | | | | | |
|---|---|---|---|---|---|---|---|---|---|---|---|---|---|
| | | | | | | 1 | 1 | 3 | 1 | | | | |
| | | 1 | 13 | 21 | 17 | 7 | 3 | 1 | | | | | |
| | | | | | 7 | 3 | 5 | 13 | 33 | 83 | | | |
| | | | | | | 17 | 7 | 9 | 23 | 29 | 73 | 183 | …∞ |
| …∞ | 651 | 521 | 167 | 67 | 27 | 11 | 9 | 23 | | …∞ | | | |
| | | | | | …∞ | 27 | 11 | 7 | 9 | …∞ | | | |
| | | | | | | 1 | 13 | 33 | | | | | |
| 69 | 221 | 177 | 71 | 57 | 23 | 37 | 15 | 19 | 3 | | | | |
| | | | | | | 21 | 17 | 43 | 27 | 17 | | | |
| 1417 | 567 | 227 | 91 | 73 | 117 | 47 | 19 | 3 | | | | | |
| | | 172917 | 69167 | 27667 | 11067 | 4427 | | | | | | | |
| | | | | | 1 | 13 | 21 | 53 | 133 | 333 | 833 | 2083 | …∞ |
| | | | | | …∞ | 57 | 23 | 29 | …∞ | | | | |
| | | | | 3 | 77 | 31 | 25 | 63 | 79 | 99 | 31 | 39 | …∞ |
| | | …∞ | 521 | 417 | 167 | 67 | 27 | 17 | | | | | |
| | | | | …∞ | 11 | 9 | 29 | 73 | …∞ | | | | |
| | | | | | | 77 | 31 | 39 | …∞ | | | | |
| …∞ | 617 | 99 | 317 | 127 | 51 | 41 | 33 | 83 | | | | | |
| …∞ | 1057 | 423 | 677 | 271 | 217 | 87 | 35 | 11 | 7 | …∞ | | | |
| | | | | | …∞ | 23 | 37 | 93 | 233 | 583 | 1823 | 2279 | …∞ |
| …∞ | 1177 | 471 | 377 | 151 | 121 | 97 | 39 | 39 | 49 | 123 | 77 | 483 | …∞ |
| | | | | | …∞ | 51 | 41 | …∞ | | | | | |
| 13021 | 10417 | 4167 | 1667 | 667 | 267 | 107 | 43 | 27 | 17 | 43 | | | |
| | | | | | | 7 | 45 | …∞ | | | | | |
| | | | | | | 117 | 47 | …∞ | | | | | |
| | | | | | 19 | 61 | 49 | …∞ | | | | | |
| | | | | | …∞ | 127 | 51 | 1 | | | | | |
| | | | | | …∞ | 33 | 53 | …∞ | | | | | |
| | …∞ | 2667 | 1067 | 427 | 171 | 137 | 55 | …∞ | | | | | |
| | | | | …∞ | 177 | 71 | 57 | …∞ | | | | | |
| 1117 | 447 | 179 | 573 | 917 | 367 | 147 | 59 | …∞ | | | | | |
| | | | | | 17 | 109 | | | | | | | |
| | | | | | | 19 | 61 | …∞ | | | | | |
| | | | | | 49 | 157 | 63 | $\cdots\infty$ | | | | | |
| | | | | 63 | 101 | 81 | 65 | 163 | 51 | 1 | | | |
| | | | | | $\cdots\infty$ | 167 | 67 | $\cdots\infty$ | | | | | |
| | | | | | $\cdots\infty$ | 43 | 69 | $\cdots\infty$ | | | | | |
| | | | | | $\cdots\infty$ | 177 | 71 | $\cdots\infty$ | | | | | |
| | | | | | | 91 | 73 | $\cdots\infty$ | | | | | |
| 89 | 4557 | 1823 | 2917 | 1167 | 467 | 187 | 75 | $\cdots\infty$ | | | | | |
| | | | $\cdots\infty$ | 27 | 173 | 277 | | | | | | | |
| | | | | | | 3 | 77 | $\cdots\infty$ | | | | | |
| $\cdots\infty$ | 27 | …. | 767 | 307 | 123 | 197 | 79 | $\cdots\infty$ | | | | | |
| | | | | | 63 | 101 | 81 | $\cdots\infty$ | | | | | |
| | | $\cdots\infty$ | 807 | 323 | 517 | 207 | 83 | 13 | | | | | |
| | | | | | $\cdots\infty$ | 53 | 85 | $\cdots\infty$ | | | | | |
| | | | | | $\cdots\infty$ | 217 | 87 | $\cdots\infty$ | | | | | |
| | | | | | $\cdots\infty$ | 111 | 89 | $\cdots\infty$ | | | | | |
| | | | | | | 227 | 91 | $\cdots\infty$ | | | | | |
| | | | | $\cdots\infty$ | 9 | 29 | 93 | $\cdots\infty$ | | | | | |
| | | | | $\cdots\infty$ | 37 | 237 | 95 | $\cdots\infty$ | | | | | |
| | | | | | $\cdots\infty$ | 121 | 97 | 243 | 19 | 3 | 1 | | |
| | | $\cdots\infty$ | 27 | … | 617 | 247 | 99 | $\cdots\infty$ | | | | | |

**Figure 3.2**

The regularities of transformations (3.7) are summarized in the form of a diagram in Figure 3.2. We see that the conjecture $C_{5q+1}$ has three periodic completion cycles $cycle_{1\leftrightarrow 1}$, $cycle_{13\leftrightarrow 13}$, $cycle_{17\leftrightarrow 17}$, similarly to the conjecture $C_{3q-1}$ [2,6]. One of the cycles of completion of the process $C_{5q+1}$ with the minimum value of an odd number $q=13$ is shown in Figure 3.1b. The conjecture $C_{5q-1}$ has only a periodic cycle $cycle_{1\leftrightarrow 1}$.

Let's check the correctness of formula (27) [1] on the example of a periodic cycle $cycle_{13\leftrightarrow 13}$ (Figure 3.1b). For the trajectory of the transformation of the number $q=13$ in the $n\to\infty$ direction we have: $q_{odd,\,initial}=q_{odd,5.1.1}=13, k_{q/2}=5$ , $m_{q/2}=1, s_{q/2}=1$.

Then from formula (1.11) we obtain that $q_{odd5,1,1}=\dfrac{\dfrac{13\cdot 2^5-1}{5}2^1-1}{5}2^1-1=((13\cdot(2^5-1)2-1)2-1=13,$

that is, is equal to the initial odd number $q_{initial}=13$ in the $cycle_{13\leftrightarrow 13}$.

Let's check the minimum values of odd numbers in cycles with another method $cycle_{1\leftrightarrow 1}$, $cycle_{13\leftrightarrow 13}$, $cycle_{17\leftrightarrow 17}$:

$q_{odd,end} =: (((5q_{odd,end}+1):2)\cdot 5+1):2:2:2:2 \Rightarrow\ 32q_{odd,end}-2=25q_{odd,end}+5 \Rightarrow q_{odd,end}=1$.

$$q_{odd,end}=((((((5q_{odd,end}+1):2)\cdot 5+1):2)5+1):2)5+1):2:2:2:2:2\Rightarrow 32q_{odd,end}=((((((5q_{odd,end}+1)$$

$$:2)\cdot 5+1):2)5+1):2)5+1128q_{odd,end}-14=(5n+1)\cdot 25\Rightarrow 3q_{odd,end}=39\Rightarrow q_{odd,end}=13$$

$$128q_{odd,end}-14=(5q_{odd,end}+1)\cdot 25\Rightarrow 3q_{odd,end}=39\Rightarrow q_{odd,end}=13,$$

$$q_{odd,end}=:((((((5q_{odd,end}+1):2)\cdot 5+1):2:2:2)5+1):2:2:2\Rightarrow 64q_{odd,end}-13=((5n+1):2)\cdot 25\Rightarrow$$

$$128q_{odd,end}-26=(5q_{odd,end}+1)\cdot 25\Rightarrow 3q_{odd,end}=51\Rightarrow q_{odd,end}=17.$$

In addition, in contrast to conjectures $C_{1,3q\pm1}$, in conjectures $C_{5q\pm1}$, growing sequences (indicated by the symbol $\cdots\infty$ in Figure 3.2) appear within the experiment of 300,000 iterations. Here, odd numbers are marked in red, which have already occurred in sequences of a smaller odd number $q$.

**Conclusion:** The regularities of transformations of recurrent numbers $K^{\pm}_{5,\theta,n}$ form the numbers $m(p)_{5,\theta,n}$ nodes on the sequences $\theta\cdot 2^n$ with a period $T_{a=5}=2^2$ and the rules of their transformations, from which the Jacobsthal tree branches out, forming $CT_{5q\pm1}$ in the reverse $n\to 0$ direction. The conjecture $C_{5q\pm1}$ will wind up with cycles $cycle_{1\leftrightarrow1}$ $(C_{5q\pm1})$ and $cycle_{13,17\leftrightarrow13,17}$ $(C_{5q+1})$ isolated from each other and infinitely growing trajectories. In the first track

$$one\ \ track:\ \ q_{odd}\to 5q_{odd}\pm1\to(5q_{odd}\pm1)/2^k=q^{\pm}_{odd,new}$$

of the $CT_{5q\pm1}$, inquality of type

$$q_{odd}\langle q^{\pm}_{odd,new}\quad \text{and} \quad q_{odd}\rangle q^{\pm}_{odd,new}$$

are fulfilled. Numbers $q^{+}_{odd,new}$

| $q_{odd}$ | 1 | 3 | 5 | 7 | 9 | 11 | 13 | 15 | 17 | 19 | 21 | 23 | 25 | … |
|---|---|---|---|---|---|---|---|---|---|---|---|---|---|---|
| | ↓ | ↓ | ↓ | ↓ | ↓ | ↓ | ↓ | ↓ | ↓ | ↓ | ↓ | ↓ | ↓ | |
| $q^{+}_{odd,new}$ | 3 | 1 | 13 | 9 | 23 | 7 | 33 | 19 | 43 | 3 | 53 | 29 | 63 | |

are divided into four groups.

**Case 4.** $a=7$. In this model, formula (5) is written in the form:

$$K^{\pm}_{7,\theta,n}=\frac{1}{7}\left[\theta\cdot 2^n\pm(-1)^n\right]. \qquad (4.1)$$

As follows from calculations (4.2)

| θ_i \2^n | **$2^0$** | **$2^1$** | **$2^2$** | **$2^3$** | **$2^4$** | **$2^5$** | **$2^6$** | **$2^7$** | **$2^8$** | | OEIS [26] |
|---|---|---|---|---|---|---|---|---|---|---|---|
| 1 | 0 | - | - | 1 | | - | 9 | - | - | $K^{\pm}_{7,1,n}$ | *A06752* |
| 3 | - | 1 | - | - | 7 | | - | 55 | - | $K^{\pm}_{7,5,n}$ | *A070997* |
| 5 | - | - | 3 | - | - | 23 | - | - | 183 | $K^{\pm}_{7,5,n}$ | *Unkn* |
| 9 | | - | 5 | - | - | 41 | - | - | 329 | $K^{\pm}_{7,9,n,}$ | *Unkn* |
| 11 | - | 3 | - | - | 25 | - | - | 201 | - | $K^{\pm}_{7,11,n}$ | *A289164* |

(4.2)

for the numbers $K^{\pm}_{7,\theta,n}$ the recurrence relation holds

$$K^{\pm}_{7,\theta,n+1} = 8K^{\pm}_{7,\theta,n} \pm 1 \ . \qquad (4.3)$$

The analysis of numbers divided by the bisection by

$$K^{-}_{7,\theta,n} = \frac{1}{7}\left[\theta \cdot 2^n - (-1)^n\right](a) \ \text{ and } \ K^{+}_{7,\theta,m} = \frac{1}{7}\left[\theta \cdot 2^n + (-1)^n\right](b), \qquad (4.4)$$

signs (4.5)

| $\theta_i \backslash 2^n$ | $2^0$ | $2^1$ | $2^2$ | $2^3$ | $2^4$ | $2^5$ | $2^6$ | $2^7$ | $2^8$ | $2^9$ | |
|---|---|---|---|---|---|---|---|---|---|---|---|
| 11 | - | - | - | - | 25 | - | - | - | - | - | $K^{+}_{7,11,n}$ |
| 9 | | - | 5 | - | | - | - | - | 329 | - | $K^{+}_{7,9,n}$ |
| 5 | - | | - | - | - | 23 | - | 1463 | - | - | $K^{+}_{7,5,n}$ |
| 3 | - | 1 | - | - | - | - | - | 55 | - | - | $K^{+}_{7,3,n}$ |
| 1 | 0 | - | - | - | - | - | 9 | - | | | $K^{+}_{7,1,n}$ |
| $\theta_i \backslash 2^n$ | $2^0$ | $2^1$ | $2^2$ | $2^3$ | $2^4$ | $2^5$ | $2^6$ | $2^7$ | $2^8$ | $2^9$ | |
| 1 | - | - | - | 1 | - | - | | - | - | 73 | $K^{-}_{7,1,n}$ |
| 3 | - | - | - | - | 7 | - | - | - | - | - | $K^{-}_{7,3,n}$ |
| 5 | - | - | 3 | | - | - | - | - | 183 | | $K^{-}_{7,5,n}$ |
| 9 | - | - | | - | - | 41 | - | - | - | - | $K^{-}_{7,9,n,}$ |
| 11 | - | 3 | - | - | - | - | - | -201 | - | - | $K^{-}_{7,11,n}$ |

(4.5)

shows that the recurring formulas

$$K^{\pm}_{7,\theta,m+4} = 64K^{\pm}_{7,\theta,m} \pm 9 \qquad (4.6)$$

are fulfilled for the numbers in the rows.

Taking (9) into account, we write down the formulas for the numbers $m(p)_{7,\theta,n}$:

$$m_{7,\theta,n} = \frac{1}{7}\left[\theta \cdot 2^n - 1\right] \ \ (a) \quad \text{and} \quad p_{7,\theta,n} = \frac{1}{7}\left[\theta \cdot 2^n + 1\right] \ \ (b), \qquad (4.7)$$

and calculate the numbers for $\theta = 1 \div 29$:

| $\theta_i \backslash 2^n$ | $2^0$ | | $2^1$ | $2^2$ | $2^3$ | $2^4$ | $2^5$ | $2^6$ | $2^7$ | $2^8$ | $2^9$ | |
|---|---|---|---|---|---|---|---|---|---|---|---|---|
| 29 | - | | - | - | - | - | - | - | - | - | - | $p_{7,29,n}$ |
| 27 | 4 | | | | 31 | | | 247 | | | 1975 | $p_{7,27,n}$ |
| 25 | - | | - | - | - | - | - | - | - | - | - | $p_{7,23,n}$ |
| 23 | - | | - | - | - | - | - | - | - | - | - | $p_{7,23,n}$ |
| 19 | | | | 11 | | | 87 | | | 695 | | $p_{7,19,n}$ |
| 17 | - | | 5 | - | - | 39 | - | - | 341 | - | - | $p_{7,17,n}$ |
| 15 | - | | - | - | - | - | - | - | - | - | - | $p_{7,15,n}$ |
| 13 | 2 | | | | 15 | | | 119 | | | 951 | $p_{7,13,n}$ |
| 11 | - | | - | - | - | - | - | - | - | - | - | $p_{7,11,n}$ |
| 9 | - | | - | - | - | - | - | - | - | - | - | $p_{7,9,n}$ |
| 5 | | | | 3 | | | 23 | | | 183 | | $p_{7,5,n}$ |
| 3 | | | 1, k=1 | | | 7, k=2 | | | 55, k=3 | | | $p_{7,3,n}$ |
| 1 | - | | - | - | - | - | - | - | - | - | - | $p_{7,1,n}$ |
| $\theta_i \backslash 2^n$ | $2^0$ | | $2^1$ | $2^2$ | $2^3$ | $2^4$ | $2^5$ | $2^6$ | $2^7$ | $2^8$ | $2^9$ | |
| 1 | [0],k=1 | | | | 1, k=2 | | | 9, k=3 | | | 73, k=3 | $m_{7,1,}$ |

| | | | | | | | | | | | | $n$ |
|---|---|---|---|---|---|---|---|---|---|---|---|---|
| 3 | - | | - | - | - | - | - | - | - | - | - | $m_{7,3,n}$ |
| 5 | - | | - | - | - | - | - | - | - | - | - | $m_{7,5,n}$ |
| 9 | | | | 5 | | | 41 | | | [329] | | $m_{7,9,n}$ |
| 11 | | | 3 | | | 25 | | | 201 | | | $m_{7,11,n}$ |
| 13 | - | | - | - | - | - | - | - | - | - | - | $m_{7,13,n}$ |
| 15 | 2 | | | | 17 | | | 137 | | | 1097 | $m_{7,15,n}$ |
| 17 | - | | - | - | - | - | - | - | - | - | - | $m_{7,17,n}$ |
| 19 | - | | - | - | - | - | - | - | - | - | - | $m_{7,97,n}$ |
| 23 | | | | 13 | | | 105 | | | 841 | | $m_{7,23,n}$ |
| 25 | | | 7 | | | 57 | | | 457 | | | $m_{7,25,n}$ |
| 27 | - | | - | - | - | - | - | - | - | - | - | $m_{7,27,n}$ |
| 29 | 4 | | | | 33 | | | 265 | | | 2121 | $m_{7,29,n}$ |

(4.8)

Square brackets highlight numbers that are multiples of seven. For numbers $m(p)_{7,\theta,k}$ in lines, the conversion rules are followed

| | | | | | | | | | |
|---|---|---|---|---|---|---|---|---|---|
| **$1\cdot2^{16}$** | | | | | | | | | |
| **$1\cdot2^{15}$** | **$4681\cdot2^{0}$** | | | | | | | | |
| **$1\cdot2^{14}$** | | | | | | | | | |
| **$1\cdot2^{13}$** | | **$167\cdot2^{0}$** | | | **$1337\cdot2^{0}$** | | | **10697** | |
| **$1\cdot2^{12}$** | **$585\cdot2^{0}$** | **$585\cdot2^{1}$** | **$585\cdot2^{2}$** | **$585\cdot2^{3}$** | **$585\cdot2^{4}$** | **$585\cdot2^{5}$** | **$585\cdot2^{6}$** | **$585\cdot2^{7}$** | **$585\cdot2^{8}$** |
| **$1\cdot2^{11}$** | | | | | | | | | |
| **$1\cdot2^{10}$** | | | | | | | | | |
| **$1\cdot2^{9}$** | **$73\cdot2^{0}$** | **$73\cdot2^{1}$** | **$73\cdot2^{2}$** | **$73\cdot2^{3}$** | **$73\cdot2^{4}$** | **$73\cdot2^{5}$** | **$73\cdot2^{6}$** | **$73\cdot2^{7}$** | **$73\cdot2^{8}$** |
| **$1\cdot2^{8}$** | | | | | | | | | |
| **$1\cdot2^{7}$** | | | | | | | | | |
| **$1\cdot2^{6}$** | **$9\cdot2^{0}$** | **$9\cdot2^{1}$** | **$9\cdot2^{2}$** | **$9\cdot2^{3}$** | **$9\cdot2^{4}$** | **$9\cdot2^{5}$** | **$9\cdot2^{6}$** | **$9\cdot2^{7}$** | **$9\cdot2^{8}$** |
| **$1\cdot2^{5}$** | | | **$5\cdot2^{0}$** | | | **$41\cdot2^{0}$** | | | **$329\cdot2^{0}$** |
| **$1\cdot2^{4}$** | | | **$5\cdot2^{1}$** | | | **$41\cdot2^{1}$** | | | |
| **$1\cdot2^{3}$** | **$1\cdot2^{0}$** | | **$5\cdot2^{2}$** | | | **$41\cdot2^{2}$** | | | |
| **$1\cdot2^{2}$** | | | **$5\cdot2^{3}$** | | | | | | |
| **$1\cdot2^{1}$** | | | | | | | | | |
| **$1\cdot2^{0}$** | **$0\cdot2^{0}$** | | | | | | | | |

**Figure 4.1**

$$m_{7,\theta,k+1} = 8m_{7,\theta,k} + 1 \text{ (a)} \quad \text{and} \quad p_{7\theta,k+1} = 8p_{7,\theta,k} - 1 \quad \text{(b)}, \qquad (4.9)$$

Rules (4.9) are analogous to rules of type $4x \pm 1$ in conjecture $3q \pm 1$. According (4.7), we have:

$$m(p)_{7,\theta,k+1} - m(p)_{7,\theta,k} = 8m(p)_{7,\theta,k} \pm 1 - m(p)_{5,\theta,k} = 7m(p)_{7,\theta,k} \pm 1. \qquad (4.10)$$

Formula (4.10) agrees with formulas (4.7). An illustration of the application of the numbering rule $m(p)_{7,\theta,k}$ by index $k$ is shown for the numbers $m(p)_{7,1(3),k}$.

Figure 4.1 shows the role in the formation of the JT of the regularities of the transformation of recurrent numbers according to the algorithm (4.7). Here, $\theta \cdot 2^n$ sequences that branch further are shown in yellow, and $\theta \cdot 2^n$

**7x-1** **7x+1**

| | | | | | | q | | | | | | |
|---|---|---|---|---|---|---|---|---|---|---|---|---|
| $\cdots\infty$ | 103 | 59 | 17 | 5 | 3 | 1 | 1 | | | | | |
| | | | | $\cdots\infty$ | 3 | 3 | 11 | 39 | 137 | 15 | 53 | $\cdots\infty$ |
| | | | | $\cdots\infty$ | 17 | 5 | 9 | 1 | | | | |
| | | | | $\cdots\infty$ | 3 | 7 | 25 | 11 | $\cdots\infty$ | | | |
| $\cdots\infty$ | 143 | 41 | 47 | 27 | 31 | 9 | 1 | | | | | |
| $\cdots\infty$ | 703 | 201 | 115 | 19 | 11 | 11 | 39 | $\cdots\infty$ | | | | |
| $\cdots\infty$ | 6723 | 1921 | 549 | 157 | 45 | 13 | 13 | 23 | 81 | 71 | 249 | $\cdots\infty$ |
| | | | $\cdots\infty$ | 45 | 13 | 15 | 53 | $\cdots\infty$ | | | | |
| | | | $\cdots\infty$ | 103 | 59 | 17 | 15 | $\cdots\infty$ | | | | |
| $\cdots\infty$ | 615 | 703 | 201 | 115 | 33 | 19 | 67 | 235 | 823 | 2881 | 2521 | $\cdots\infty$ |
| $\cdots\infty$ | 341 | 195 | 223 | 255 | 73 | 21 | 37 | 65 | 57 | 25 | $\cdots\infty$ | |
| | | | | $\cdots\infty$ | 5 | 23 | 81 | $\cdots\infty$ | | | | |
| | | | $\cdots\infty$ | 19 | 87 | 25 | 11 | $\cdots\infty$ | | | | |
| | | | | $\cdots\infty$ | 47 | 27 | 95 | 333 | 583 | 2041 | 893 | $\cdots\infty$ |
| $\cdots\infty$ | 7563 | 2161 | 1235 | 353 | 101 | 29 | 51 | 179 | 627 | 2195 | 7683 | $\cdots\infty$ |
| | | | | $\cdots\infty$ | 27 | 31 | 109 | 191 | 669 | 1171 | 4099 | $\cdots\infty$ |
| $\cdots\infty$ | 269 | 615 | 703 | 201 | 115 | 33 | 29 | 51 | $\cdots\infty$ | | | |
| $\cdots\infty$ | 2281 | 2607 | 745 | 213 | 61 | 35 | 123 | 431 | 1509 | 2641 | 2311 | $\cdots\infty$ |
| $\cdots\infty$ | 9663 | 2761 | 789 | 451 | 129 | 37 | 65 | $\cdots\infty$ | | | | |
| | | | | $\cdots\infty$ | 17 | 39 | 137 | 15 | $\cdots\infty$ | | | |
| | | | | $\cdots\infty$ | 143 | 41 | 9 | 1 | | | | |
| $\cdots\infty$ | 2803 | 801 | 229 | 131 | 75 | 43 | 151 | 529 | 926 | 1621 | 2837 | $\cdots\infty$ |
| $\cdots\infty$ | 11765 | 6723 | 1921 | 549 | 157 | 45 | 79 | 277 | 485 | 849 | 743 | $\cdots\infty$ |
| | | | $\cdots\infty$ | 143 | 41 | 47 | 165 | 289 | 253 | 443 | 1551 | $\cdots\infty$ |
| $\cdots\infty$ | 5603 | 1601 | 523 | 299 | 171 | 49 | 43 | 151 | 529 | 463 | 1621 | $\cdots\infty$ |
| | | $\cdots\infty$ | 17 | 311 | 89 | 51 | 179 | 627 | 2195 | 7683 | 26891 | $\cdots\infty$ |

**Figure 4.2**

sequences that do not branch further are shown in blue. Sequences $\theta \cdot 2^n$ whose parameter $\theta$ is a multiple of seven are shown in red. An illustration of the algorithm (1) for model (4.1) is shown in Figure 4.2.

Therefore, in contrast to models of recurrent numbers $K^{\pm}_{1,3,5,\theta,n}$, the numbers $p_{7,1,k}$ on the sequence $1 \cdot 2^n$ do not form nodes. Therefore, in the Collatz problem, the $C_{7q-1}$ function is not active, and the corresponding Collatz sequences are isolated from the $1 \cdot 2^n$ sequence, that is, increasing, which confirms the results of the calculations presented in the form of the left part of the diagram in Figure 4.2. For the function $C_{7q+1}$ forms a certain group of numbers 5,9,41, etc., with a single periodic cycle of $cycle_{1\leftrightarrow 1}$ complection. In addition, periodicity of groups of conjectures is highlighted in (4.8).

Let's check the formula (1.16) for the number $q_{odd,\,initial}=41$ with trajectory parameters: $k_{q/2}=6$, $m_{q/2}=5$, $q_{odd,end}=1$. According to the formula (1.16), we have: $q_{odd,end}=((1\cdot 2^6-1)2^5/7-1)/7=41$ which is equal to the initial odd number $q_{odd,initial}=41$.

**Conclusion:** The regularities of transformations of recurrent numbers $K^{\pm}_{7,\theta,n}$ form the numbers $m(p)_{7,\theta,n}$ nodes on the sequences $\theta\cdot 2^n$ with a period $T_{a=7}=3$ and the rules of their transformations, from which the *JT* branches out, forming $CT_{7q+1}$ in the reverse direction $C_{15q-1}$. The odd value of the period $T_{a=7}=3$ of the numbers $m(p)_{7,\theta,n}$ forms isolated from $\theta\cdot 2^n$ sequence the $CT_{7q-1}$ and trajectories the $C_{7q+1}$ without periodic oscillations or infinite growth.

In the first track $one\ \ track$: $q_{odd}\rightarrow 7q_{odd}+1\rightarrow(7q_{odd}+1)/2^k=q^+_{odd,new}$ of the Collatz trajectory $C_{7q+1}$, inquality of type $q_{odd}\ \langle\ q^+_{odd,new}$ and $q_{odd}\ \rangle\ q^+_{odd,new}$ are fulfilled. Numbers $q^+_{odd,new}$

| $q_{odd}$ | 1 | 3 | 5 | 7 | 9 | 11 | 13 | 15 | 17 | 19 | 21 | 23 | 25 | … |
|---|---|---|---|---|---|---|---|---|---|---|---|---|---|---|
| | ↓ | ↓ | ↓ | ↓ | ↓ | ↓ | ↓ | ↓ | ↓ | ↓ | ↓ | ↓ | ↓ | |
| $q^+_{odd,new}$ | 1 | 11 | 9 | 25 | 1 | 39 | 23 | 53 | 15 | 67 | 37 | 81 | 11 | |

are divided into two groups, the Collatz trajectories of which grow infinitely, or go to unity.

**Case 5.** $a=9$. In this model, formula (5) is written in the form:

$$K^{\pm}_{9,\theta,n}=\frac{1}{9}\left[\theta\cdot 2^n\pm(-1)^n\right]. \quad (5.1)$$

Analysis of numbers divided by bisection by signs

$$K^{-}_{9,\theta,n}=\frac{1}{9}\left[\theta\cdot 2^n-(-1)^n\right](\mathrm{a})\ \ \text{and}\ \ K^{+}_{9,\theta,m}=\frac{1}{9}\left[\theta\cdot 2^n+(-1)^n\right](\mathrm{b}), \quad (5.2)$$

and calculated in (5.3)

| | | | | | | | | | | | |
|---|---|---|---|---|---|---|---|---|---|---|---|
| 11 | - | - | - | - | - | - | - | - | - | - | $K^+_{9,11,n}$ |
| 7 | | - | 3 | - | | 25 | - | - | 199 | - | $K^+_{9,7,n}$ |
| 5 | - | | - | - | - | - | - | - | - | - | $K^+_{9,5,n}$ |
| 3 | - | - | - | - | - | - | - | - | - | - | $K^+_{9,3,n}$ |
| 1 | 0 | - | - | 1 | - | - | 7 | - | - | **57** | $K^+_{9,1,n}$ |
| $\theta_i \backslash 2^n$ | $2^0$ | $2^1$ | $2^2$ | $2^3$ | $2^4$ | $2^5$ | $2^6$ | $2^7$ | $2^8$ | $2^9$ | |
| 1 | - | - | - | - | - | - | - | - | - | - | $K^-_{9,1,n}$ |
| 3 | - | - | - | - | - | - | - | - | - | - | $K^-_{9,3,n}$ |
| 5 | - | 1 | - | - | 9 | - | - | 71 | - | - | $K^-_{9,5,n}$ |
| 9 | - | - | | - | - | - | - | - | - | - | $K^-_{9,7,n,}$ |
| 11 | - | - | 5 | - | - | 39 | - | - | 313 | - | $K^-_{9,11,n}$ |

(5.3)

shows that recursive formulas are executed for the numbers in the rows

$$K^{\pm}_{9,\theta,m+3} = 8K^{\pm}_{9,\theta,m} + \begin{cases} -1 & \text{if } m \text{ is odd}, \\ +1 & \text{if } m \text{ is even}. \end{cases} \qquad (5.4)$$

Taking (9) into account, we write down the formulas:

$$m_{9,\theta,n} = \frac{1}{9}\left[\theta \cdot 2^n - 1\right] \quad \text{(a)} \qquad \text{and} \qquad p_{9,\theta,n} = \frac{1}{9}\left[\theta \cdot 2^n + 1\right] \quad \text{(b)}, \qquad (5.5)$$

for numbers $m(p)_{9,\theta,n}$ and calculate them for $\theta = 1 \div 11$:

| | | | | | | | | | | | |
|---|---|---|---|---|---|---|---|---|---|---|---|
| 11 | - | - | 5 | - | - | - | - | - | 313 | - | $p_{9,11,n}$ |
| 7 | - | - | - | - | - | 2 | - | - | - | - | $p_{9,7,n}$ |
| 5 | - | - | - | - | 9 | - | - | - | - | - | $p_{9,5,n}$ |
| 3 | - | - | - | - | - | - | - | - | - | - | $p_{9,3,n}$ |
| 1 | - | - | - | 1, k=1 | - | - | - | - | - | 57, k=2 | $p_{9,1,n}$ |
| $\theta_i$ \ $2^n$ | $2^0$ | $2^1$ | $2^2$ | $2^3$ | $2^4$ | $2^5$ | $2^6$ | $2^7$ | $2^8$ | $2^9$ | |
| 1 | 0,k=1 | | | | | | 7, k=2 | | | | $m_{9,1,n}$ |
| 3 | - | - | - | - | - | - | - | - | - | - | $m_{9,3,n}$ |
| 5 | - | 1 | - | - | - | - | - | 71 | - | - | $m_{9,5,n}$ |
| 7 | - | - | 3 | - | - | -- | - | - | 199 | - | $m_{9,7,n}$ |
| 11 | | | | | | 39 | | | | | $m_{9,11,n}$ |

(5.6)

For numbers $m(p)_{7,\theta,k}$ in lines, the conversion rules are followed

$$m_{7,\theta,k+1} = 64m_{7,\theta,k} + 7 \quad \text{(a)} \qquad \text{and} \qquad p_{7\theta,k+1} = 64p_{7,\theta,k} - 7 \quad \text{(b)}, \qquad (5.7)$$

According (5.7), we have:

$$m(p)_{9,\theta,k+1} - m(p)_{9,\theta,k} = 64m(p)_{9,\theta,k} \pm 7 - m(p)_{9,\theta,k} = 7(9m(p)_{7,\theta,k} \pm 1). \qquad (5.8)$$

Formula (5.8) agrees with formulas (5.5).

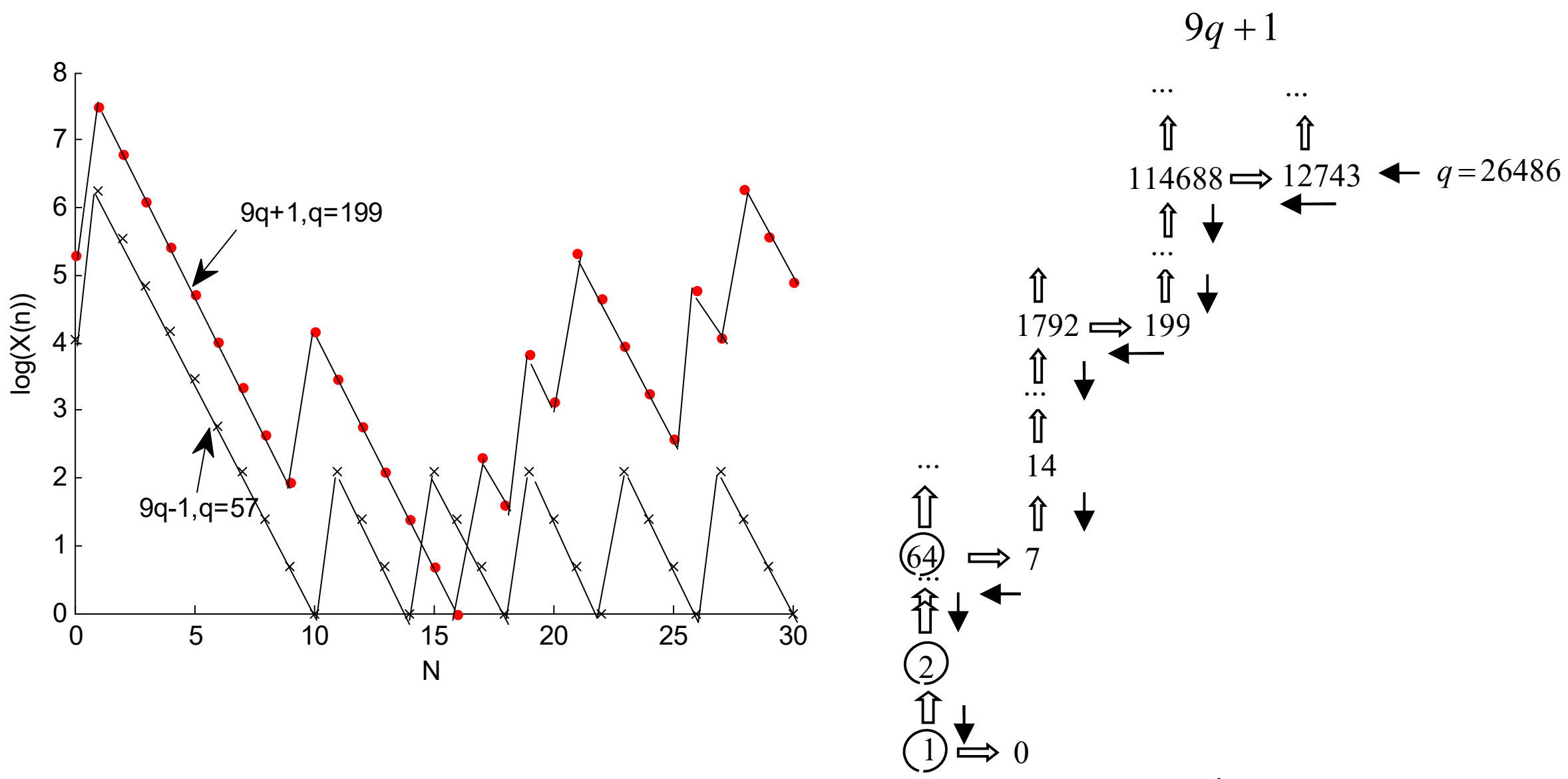

a b

**Figure 5.1**

An illustration of the application of the rule for numbering numbers $m(p)_{7,\theta,k}$ by index $k$ is shown in (5.6) for numbers $m(p)_{9,1,k}$. The rules (5.7) are analogous to the rules of the type in $4x \pm 1$ conjectures $3q \pm 1$.

However, unlike previous models of conjectures $C_{(a\le 7)q\pm 1}$ (recurring numbers $m(p)_{a\le 7,\theta,n}$), the conjecture $C_{9q+1}$ has a fundamentally new property, namely, the Collatz trajectory $C_{9q+1}$, having reached unity, unlike the conjecture $C_{9q-1}$, grows infinitely (Figure 5.1a). This figure also shows the trajectory of the Collatz conjecture of $q = 26486$.

The increasing sequence of conjecture $C_{9q+1}$ is observed for all numbers localized on the Jacobsthal tree constructed in Figure 5.2b. For others

| $p_{17,1,n}$ | $p_{15,1,n}$ | $p_{13,1,n}$ | $p_{11,1,n}$ | $p_{9,1,n}$ | $2^n$ | $m_{9,1,n}$ | $m_{11,1,n}$ | $m_{13,1,n}$ | $m_{15,1,n}$ | $p_{17,1,n}$ |
|---|---|---|---|---|---|---|---|---|---|---|
| | | | | | $2^{16}$ | | | | | 3855 |
| | | | 2979 | | $2^{15}$ | | | | | |
| | | | | | $2^{14}$ | | | | | |
| | - | | | | $2^{13}$ | | | | | |
| 241 | - | | | | $2^{12}$ | 455 | | 315 | **273** | |
| | - | | | | $2^{11}$ | | | | | |
| | - | | | | $2^{10}$ | | **93** | | | |
| | - | | | 57 | $2^{9}$ | | | | | |
| | - | | | | $2^{8}$ | | | | 17 | 15 |
| | - | | | | $2^{7}$ | | | | | |
| | - | 5 | | | $2^{6}$ | 7 | | | | |
| | - | | 3 | | $2^{5}$ | | | | | |
| 1 | - | | | | $2^{4}$ | | | | 1 | |
| | - | | | 1 | $2^{3}$ | | | | | |
| | - | | | | $2^{2}$ | | | | | |
| | - | | | | $2^{1}$ | | | | | |
| | - | | | | $2^{0}$ | 0 | 0 | 0 | **0** | 0 |
| $2^8x-15$ | | $2^{12}x-315$ | $2^{10}x-93$ | $2^6x-7$ | | $2^6x+7$ | $2^{10}x+93$ | $2^{12}x+315$ | $2^4x+1$ | $2^8x+15$ |

(5.9)

conjectures, the same patterns described are displayed in the same colors. We also see that fractional values accept numbers $p_{15,1,k}$, therefore, the trajectories of number conjecture by function $C_{15q-1}$ are also isolated from the sequence $1 \cdot 2^n$. The bottom line (5.9) contains the formulas connecting the adjacent numbers $m(p)_{a,\theta,k}$ and $m(p)_{a,\theta,k+1}$ in each column.

**Conclusion**: The regularities of transformations of recurrent numbers $K^{\pm}_{9,\theta,n}$ form the numbers $m(p)_{9,\theta,n}$ nodes on the sequences $\theta \cdot 2^n$ with a period $T_{a=7} = 2^6$ and the rules of their transformations, from which the *JT* branches out, forming infinitely growing $CT_{9q\pm 1}$.

In the first track $one\ track: \ q_{odd} \to 9q_{odd} \pm 1 \to (9q_{odd} \pm 1)/2^k = q^{\pm}_{odd,new}$ of the $CT_{9q\pm 1}$, inquality of type $q_{odd} \langle\, q^{\pm}_{odd,new}$ and $q_{odd} \,\rangle\, q^{\pm}_{odd,new}$ are fulfilled. Numbers $q^{+}_{odd,new}$

| $q_{odd}$ | 1 | 3 | 5 | 7 | 9 | 11 | 13 | 15 | 17 | 19 | 21 | 23 | 25 | … |
|---|---|---|---|---|---|---|---|---|---|---|---|---|---|---|
| | ↓ | ↓ | ↓ | ↓ | ↓ | ↓ | ↓ | ↓ | ↓ | ↓ | ↓ | ↓ | ↓ | |
| $q^{+}_{odd,new}$ | 5 | 7 | 23 | 1 | 41 | 25 | 59 | 17 | 77 | 43 | 95 | 13 | 113 | |

are formed one group, the $CT_{9q\pm1}$ of which grow infinitely.

**General conclusions**

1.Positive integers

$$K^{\pm}_{a,\theta,n} = \frac{\theta \cdot 2^n \pm (-1)^n}{a} = Integer\,,\; n \in \mathbb{N},\; a,\theta \in \mathbb{N}_{odd}\,,\; \mathbb{N} = \mathbb{N}_{odd} + \mathbb{N}_{even}$$

recurrent, of the second order and calculated according to the formulas $K^{\pm}_{a,\theta,n+2} = K^{\pm}_{a,\theta,n+1} + 2K^{\pm}_{a,\theta,n}$.

2.The bisection of numbers $K^{\pm}_{a,\theta,n}$ by signs $K^{\pm}_{a,\theta,n} = \begin{cases} K^{+}_{a,\theta,n} = \dfrac{\theta \cdot 2^n + (-1)^n}{a}, \\ K^{-}_{a,\theta,n} = \dfrac{\theta \cdot 2^n - (-1)^n}{a}, \end{cases}$ forms parameterized $a,\theta$ nodes with numbers $m(p)_{a,\theta,k} = \dfrac{\theta \cdot 2^k \pm 1}{a} = Integer \;\; k \in \mathbb{N}$. From nodes with numbers $m(p)_{a,\theta,k}$ branch $\theta \cdot 2^k$ sequences and form a *JT*, in the reverse $n \to 0$ direction of which the $CT_{aq\pm1}$ of numbers $q \in \mathbb{N}$ according to algorithm (2). The numbers $m_{a,1,n} = \dfrac{1}{a}\left[1 \cdot 2^n - 1\right]$ and $p_{a,1,n} = \dfrac{1}{a}\left[1 \cdot 2^n + 1\right]$ of nodes $1 \cdot 2^n$ sequence are equal:

| 9q-1 | 7q-1 | 5q-1 | 3q-1 | 1q-1 | $2^s$ | 1q+1 | 3q+1 | 5q+1 | 7q+1 | **9q+1** |
|---|---|---|---|---|---|---|---|---|---|---|
| | - | | 683 | 2049 | $2^{11}$ | 2045 | | | | |
| | - | 205 | | *1025* | $\mathbf{2^{10}}$ | 1023 | 341 | | | |
| 57 | - | | 171 | 513 | $2^9$ | 511 | | | 73 | |
| | - | | | *257* | $\mathbf{2^8}$ | 255 | 85 | 51 | | |
| | -- | | 43 | 129 | $\mathbf{2^7}$ | 127 | | | | |
| | - | 13 | | ***65*** | $\mathbf{2^6}$ | 63 | 21 | | 9 | 7 |
| | - | | 11 | **33** | $2^5$ | 31 | | | | |
| | - | | | *17* | $\mathbf{2^4}$ | 15 | 5 | 3 | | |
| 1 | **-** | | 3 | 9 | $2^3$ | 7 | | | 1 | |
| | - | 1 | | *5* | $\mathbf{2^2}$ | 3 | **1** | | | |
| | - | | 1 | 3 | $2^1$ | 1 | | | | |
| | - | | | *2* | $\mathbf{2^0}$ | 0 | 0 | 0 | **0** | 0 |
| $p_{9,1,n}$ | $p_{7,1,n}$ | $p_{5,1,n}$ | $p_{3,1,n}$ | $p_{,1,n}$ | | $m_{1,1,n}$ | $m_{3,1,n}$ | $m_{5,1,n}$ | $m_{7,1,n}$ | $m_{9,1,n}$ |
| $2^6x-7$ | | $2^4x-3$ | $2^2x-1$ | $2^1x-1$ | | $2^1x+1$ | $2^2x+1$ | $2^4x+3$ | $2^3x+1$ | $2^6x+7$ |

and

| 17q-1 | 15q-1 | 13q-1 | 11q-1 | 9q-1 | $\mathbf{2^n}$ | 9q+1 | 11q+1 | 13q+1 | 15q+1 | 17q+1 |
|---|---|---|---|---|---|---|---|---|---|---|
| | | | | | $\mathbf{2^{16}}$ | | | | | 3855 |
| | | | 2979 | | $\mathbf{2^{15}}$ | | | | | |
| | | | | | $\mathbf{2^{14}}$ | | | | | |
| | - | | | | $\mathbf{2^{13}}$ | | | | | |
| 241 | - | | | | $\mathbf{2^{12}}$ | 455 | | 315 | **273** | |
| | - | | | | $\mathbf{2^{11}}$ | | | | | |
| | - | | | | $\mathbf{2^{10}}$ | | **93** | | | |
| | - | | | 57 | $\mathbf{2^9}$ | | | | | |
| | - | | | | $\mathbf{2^8}$ | | | | 17 | 15 |
| | - | | | | $\mathbf{2^7}$ | | | | | |
| | - | 5 | | | $\mathbf{2^6}$ | 7 | | | | |
| | - | | 3 | | $\mathbf{2^5}$ | | | | | |
| 1 | **-** | | | | $\mathbf{2^4}$ | | | | 1 | |
| | - | | | 1 | $\mathbf{2^3}$ | | | | | |
| | - | | | | $\mathbf{2^2}$ | | | | | |
| | - | | | | $\mathbf{2^1}$ | | | | | |
| | - | | | | $\mathbf{2^0}$ | 0 | 0 | 0 | **0** | 0 |
| $p_{17,1,n}$ | $p_{15,1,n}$ | $p_{13,1,n}$ | $p_{11,1,n}$ | $p_{9,1,n}$ | | $m_{9,1,n}$ | $m_{11,1,n}$ | $m_{13,1,n}$ | $m_{15,1,n}$ | $m_{17,1,n}$ |
| $2^8x-15$ | | $2^{12}x+315$ | $2^{10}x-93$ | $2^6x-7$ | | $2^6x+7$ | $2^{10}x+93$ | $2^{12}x+315$ | $\mathbf{2^4x\ +1}$ | $2^8x+15$ |

The

period $T_a$ between two nodes is equal:

| $a$ | 1 | 3 | 5 | 7 | 9 | 11 | 13 | 15 | 17 |
|---|---|---|---|---|---|---|---|---|---|
| $T_a$ | 1 | 4 | 16 | 8 | 64 | 1024 | 4096 | 16 | 256 |

2. When $a \geq 5$, infinitely increasing Collatz sequences are formed. Conjectures of the unit by the Collatz function $aq+-1$ are completed by periodic cycles

$a=1$: $1 \cdot 1 - 1 \to 0,$

$$\begin{rrcases} a=1: & 1\cdot 1+1 \to 2 \to 1, \\ a=3: & 3\cdot 1-1 \to 2 \to 1, \end{rrcases} \Rightarrow \begin{cases} 1 = 2^1 - 1, \\ 3 = 2^1 + 1, \end{cases}$$

$$\left.\begin{array}{ll} a=3: & 3\cdot 1+1 \to 4 \to 2 \to 1, \\ a=5: & 5\cdot 1-1 \to 4 \to 2 \to 1, \end{array}\right\} \Rightarrow \begin{cases} 3 = 2^2 - 1, \\ 5 = 2^2 + 1, \end{cases}$$

$$\left.\begin{array}{ll} a=5: & 5\cdot 1+1 \to 6 \to 3 \to 16 \to 8 \to 4 \to 2 \to 1, \\ a=7: & 7\cdot 1-1 \to 6 \to 3 \to 20 \to 5 \to 34 \to 17 \to \ldots \to \infty, \end{array}\right\}$$

$$\left.\begin{array}{ll} a=7: & 7\cdot 1+1 \to 8 \to 4 \to 2 \to 1, \\ a=9: & 9\cdot 1-1 \to 8 \to 4 \to 2 \to 1, \end{array}\right\} \Rightarrow \begin{cases} 7 = 2^3 - 1, \\ 9 = 2^3 + 1, \end{cases}$$

$$\left.\begin{array}{ll} a=9: & 9\cdot 1+1 \to 10 \to 5 \to 46 \to 23 \to 208 \to \ldots \to \infty, \\ a=11: & 11\cdot 1-1 \to 10 \to 5 \to 54 \to 27 \to 296 \to \ldots \to \infty, \end{array}\right\}$$

$$\left.\begin{array}{ll} a=11: & 11\cdot 1+1 \to 12 \to 6 \to 3 \to 34 \to 17 \to \ldots \to \infty, \\ a=13: & 13\cdot 1-1 \to 12 \to 6 \to 3 \to 38 \to 19 \to \ldots \to \infty, \end{array}\right\}$$

$$\left.\begin{array}{ll} a=13: & 13\cdot 1+1 \to 12 \to 6 \to 3 \to 40 \to 20 \to \ldots \to \infty, \\ a=15: & 15\cdot 1-1 \to 14 \to 7 \to 85 \to 1274 \to 637 \to \ldots \to \infty, \end{array}\right\}$$

$$\left.\begin{array}{ll} a=15: & 15\cdot 1+1 \to 16 \to 8 \to 4 \to 2 \to 1, \\ a=17: & 17\cdot 1-1 \to 16 \to 8 \to 4 \to 2 \to 1, \end{array}\right\} \Rightarrow \begin{cases} 15 = 2^4 - 1, \\ 17 = 2^4 + 1, \end{cases}$$

.......... ...

if for adjacent values of $a, a+2$, the $a = 2^n \mp 1$ condition is satisfied, otherwise the Collatz sequences are increasing. Indeed, if the track of a periodic cycle has one step of $k$ iterations of $q/2$, then $q = \pm\frac{1}{2^k - a}$, whence $1 = \pm\frac{1}{2^k - a}$, if $a = 2^k \pm 1$. If $a = 2^k - 1$ $(aq+1)$, then

| $k$ | 1 | 2 | 3 | 4 | 5 | 6 | 7 | … |
|---|---|---|---|---|---|---|---|---|
| $a$ | 1 | 3 | 7 | 15 | 31 | 63 | 127 | |

If $a = 2^k + 1$ $(aq-1)$, then

| $k$ | 1 | 2 | 3 | 4 | 5 | 6 | 7 | … |
|---|---|---|---|---|---|---|---|---|
| $a$ | 3 | 5 | 9 | 17 | 33 | 65 | 129 | |

3. In $C_{a\cdot q\pm1} = \ \ if \quad q\equiv 0 \mod 2 \quad then \quad \frac{q}{2} \quad else \quad a\cdot q \pm 1$ conjectures $q\in\mathbb{N}$ numbers isolated $cycle^{3q-1}_{5(17)\leftrightarrow5(17)}$ and $cycle^{5q+1}_{13(17)\leftrightarrow13(17)}$ from the sequence $1\cdot 2^n$ cycles are formed only for two $CT_{3q-1}$ and $CT_{5q+1}$ conjectures, in which the minimum values of odd numbers are localized in relatively small values. Distribution of the numbers $q\in\mathbb{N}$ by cycles $cycle^{3q-1}_{1,5,17\leftrightarrow1,5,17}$ has a trend towards a uniform distribution. The minimum odd numbers in cycles $cycle^{3q-1}_{1\leftrightarrow1}$, $cycle^{3q-1}_{5(7)\leftrightarrow5(7)}$ and $cycle^{3q-1}_{17\leftrightarrow17}$ are localized in the interval of small values. Therefore, for the conjecture $C_{3q+1}$ the absence of similar cycles with values that are not equal to one is a "strong" argument in favor of the Collatz hypothesis.

## Appendix 1.

Investigation of the problem $aq\pm1$ indicates that the regularities of transforming natural numbers $q\in\mathbb{N}$ in the $n\to0$ direction are correlated with the regularities of periods in $T_a$ of Jacobsthal numbers $m(p)_{a,\theta,n}$. Let's demonstrate that the period $T_a$ correlates with the maximum number $k$ of doublings of the number $q$ between two consecutive transitions

$$\frac{q\cdot 2^n \pm 1}{a} = Integer \qquad \text{(A1-1)}$$

in the $n\to\infty$ direction. To achieve this, for the conjectures $C_{3,5,7\pm1}$, we will construct the so-called Jacobsthal value matrices, as shown in (A1-1) - (A1-3) for odd numbers forming nodes of the $JT$.

**Jacobsthal matrices**

3x+1 | 3x-1

| | 1 | 3 | 5 | 7 | 9 | | 1 | 3 | 5 | 7 | 9 |
|---|---|---|---|---|---|---|---|---|---|---|---|
| 00 | 2 | ∞ | 1 | 2 | ∞ | 00 | 1 | ∞ | 2 | 1 | ∞ |
| 10 | 1 | 2 | ∞ | 1 | 2 | 10 | 2 | 1 | ∞ | 2 | 1 |
| 20 | ∞ | 1 | 2 | ∞ | 1 | 20 | ∞ | 2 | 1 | ∞ | 2 |
| 30 | 2 | ∞ | 1 | 2 | ∞ | 30 | 1 | ∞ | 2 | 1 | ∞ |
| 40 | 1 | 2 | ∞ | 1 | 2 | 40 | 2 | 1 | ∞ | 2 | 1 |
| 50 | ∞ | 1 | 2 | ∞ | 1 | 50 | ∞ | 2 | 1 | ∞ | 2 |
| 60 | 2 | ∞ | 1 | 2 | ∞ | 60 | 1 | ∞ | 2 | 1 | ∞ |
| 70 | 1 | 2 | ∞ | 1 | 2 | 70 | 2 | 1 | ∞ | 2 | 1 |
| 80 | ∞ | 1 | 2 | ∞ | 1 | 80 | ∞ | 2 | 1 | ∞ | 2 |
| 90 | 2 | ∞ | 1 | 2 | ∞ | 90 | 1 | ∞ | 2 | 1 | ∞ |
| 100 | 1 | 2 | ∞ | 1 | 2 | 100 | 2 | 1 | ∞ | 2 | 1 |
| 110 | ∞ | 1 | 2 | ∞ | 1 | 110 | ∞ | 2 | 1 | ∞ | 2 |
| 120 | 2 | ∞ | 1 | 2 | ∞ | 120 | 1 | ∞ | 2 | 1 | ∞ |
| 130 | 1 | 2 | ∞ | 1 | 2 | 130 | 2 | 1 | ∞ | 2 | 1 |
| 140 | ∞ | 1 | 2 | ∞ | 1 | 140 | ∞ | 2 | 1 | ∞ | 2 |
| 150 | 2 | ∞ | 1 | 2 | ∞ | 150 | 1 | ∞ | 2 | 1 | ∞ |
| 160 | 1 | 2 | ∞ | 1 | 2 | 160 | 2 | 1 | ∞ | 2 | 1 |
| 170 | ∞ | 1 | 2 | ∞ | 1 | 170 | ∞ | 2 | 1 | ∞ | 2 |
| **180** | **2** | ∞ | **1** | **2** | ∞ | **180** | **1** | ∞ | **2** | **1** | ∞ |
| **…** | | | | | | **…** | | | | | |

(A1-2)

For example, if $q_{odd} = 177$ (a row with the symbol **170** and a column with the symbol **7**), then for $C_{3q+1}$ the relation $\frac{177 \cdot 2^n - 1}{3} \neq Integer$, holds, as the number $177/3 = 59$ is a multiple of three. Therefore, the sequence with the parameter $\theta = 177$ does not branch, and in (A1-2) $k_{q/2} \to \infty$. Let $q_{odd} = 121$ (arow with the symbol **120** and a column with the symbol **1**). The number $121/3 \neq Integer$ is not a multiple of three, so for $C_{3q+1}$ the first branching node of the sequence $\theta = 121$ is formed after the second doubling iteration $\frac{(2 \cdot 121) \cdot 2 - 1}{3} = 161 = Integer$.

The value of $k_{q/2} = 2$ is the maximum value of doublings of an odd number in the $n \to \infty$ direction and is equal to the period $T_3 = 2$ of the numbers $m(p)_{3,\theta,n}$. From (A1-2), it also follows that for the values of $k_{q/2}$, the diagrams of both $C_{3q\pm1}$ conjectures form diagonal-type isodirections of

$$k_{q/2} = const, \qquad \text{(A1-3)}$$

which change symmetrically $2 \leftrightarrow 1$ with the transition $C_{3q+1} \leftrightarrow C_{3q-1}$. Thus, the branching regularities of the *JT* change with the conjectures $C_{3q+1} \leftrightarrow C_{3q-1}$. Indeed, in the $C_{3q-1}$ model, in the $n \to \infty$ direction, different trees are formed, isolated from the sequence $1 \cdot 2^n$ with the beginnings $5 \cdot 2^n$ and $17 \cdot 2^n$.

In (A1-4), similar to (A1-2), Jacobsthal diagrams for type $C_{5q\pm1}$ conjectures are presented:

| 5x+1 | | | | | | 5x-1 | | | | | |
|---|---|---|---|---|---|---|---|---|---|---|---|
| | 1 | 3 | 5 | 7 | 9 | | 1 | 3 | 5 | 7 | 9 |
| **00** | **4** | **1** | ∞ | **3** | **2** | **00** | **2** | **3** | ∞ | **1** | **4** |
| **10** | **4** | **1** | ∞ | **3** | **2** | **10** | **2** | **3** | ∞ | **1** | **4** |
| **20** | **4** | **1** | ∞ | **3** | **2** | **20** | **2** | **3** | ∞ | **1** | **4** |
| **30** | **4** | **1** | ∞ | **3** | **2** | **30** | **2** | **3** | ∞ | **1** | **4** |
| **40** | **4** | **1** | ∞ | **3** | **2** | **40** | **2** | **3** | ∞ | **1** | **4** |
| **50** | **4** | **1** | ∞ | **3** | **2** | **50** | **2** | **3** | ∞ | **1** | **4** |
| **60** | **4** | **1** | ∞ | **3** | **2** | **60** | **2** | **3** | ∞ | **1** | **4** |
| **70** | **4** | **1** | ∞ | **3** | **2** | **70** | **2** | **3** | ∞ | **1** | **4** |
| **80** | **4** | **1** | ∞ | **3** | **2** | **80** | **2** | **3** | ∞ | **1** | **4** |
| **90** | **4** | **1** | ∞ | **3** | **2** | **90** | **2** | **3** | ∞ | **1** | **4** |
| **100** | **4** | **1** | ∞ | **3** | **2** | **100** | **2** | **3** | ∞ | **1** | **4** |
| **110** | **4** | **1** | ∞ | **3** | **2** | **110** | **2** | **3** | ∞ | **1** | **4** |
| **120** | **4** | **1** | ∞ | **3** | **2** | **120** | **2** | **3** | ∞ | **1** | **4** |
| **130** | **4** | **1** | ∞ | **3** | **2** | **130** | **2** | **3** | ∞ | **1** | **4** |
| **140** | **4** | **1** | ∞ | **3** | **2** | **140** | **2** | **3** | ∞ | **1** | **4** |
| **150** | **4** | **1** | ∞ | **3** | **2** | **150** | **2** | **3** | ∞ | **1** | **4** |
| **160** | 4 | 1 | ∞ | 3 | 2 | **160** | 2 | 3 | ∞ | 1 | 4 |
| **170** | 4 | 1 | ∞ | 3 | 2 | **170** | 2 | 3 | ∞ | 1 | 4 |
| **180** | 43 | 1 | ∞ | 3 | 2 | **180** | 2 | 3 | ∞ | 1 | 4 |
| **…** | | | | | | **…** | | | | | |

(A1-4)

Unlike (A1-2), here, isodirections of type (A1-3) are formed along the columns with the maximum value of $k = 4$, which equals the period $T_5 = 4$. Therefore, the regularities of $C_{5q\pm1}$ conjectures fundamentally differ from the regularities of $C_{3q\pm1}$ conjectures. Indeed, in the $a = 5$ model, additional periodic cycles isolated from the sequence $1\cdot 2^n$ are formed for the $C_{5q+1}$ conjectures. With the transition $C_{5q+1} \leftrightarrow C_{5q-1}$, they are symmetric relative to the column with the symbol 5. In the $a = 7$ model, the period $T_7 = 3$ equals an odd number. Therefore, in the matrices

| | 7x+1 | | | | | 7x-1 | | | | |
|---|---|---|---|---|---|---|---|---|---|---|
| | **1** | **3** | **5** | **7** | **9** | | **1** | **3** | **5** | **7** | **9** |
| **00** | 3 | 1 | ∞ | ∞ | 2 | **00** | ∞ | 1 | 2 | ∞ | ∞ |
| **10** | 1 | ∞ | 3 | ∞ | ∞ | **10** | 2 | 3 | ∞ | 1 | 2 |
| **20** | ∞ | 2 | 1 | ∞ | 3 | **20** | ∞ | ∞ | ∞ | 3 | ∞ |
| **30** | ∞ | ∞ | ∞ | 2 | 1 | **30** | 1 | 2 | ∞ | ∞ | 3 |

(A-5)

of Jacobsthal $k_{q/2,\max} = 3$, there are no clearly defined isodirections (A1-3) or mirror symmetry with the transition $C_{7q+1} \leftrightarrow C_{7q-1}$.

**Appendix 2.** *Jacobsthal tree construction algorithm for 3x+1 conjecture*

Let the initial number be 1 ($\theta = 1$). We double it until (2.13) is executed: 1→2→(4-1)/3=1. We get the $cycle_{1\leftrightarrow 1}^{3q+1}$ cycle, so then we double the number 4: 4→8→(16-1)/3=5. The node with the number 5 is not a multiple of three, therefore it generates a new $\theta \cdot 2^k$ sequence, with which we build the bottom branch on the *JT*, doubling the number 5: 5→10→ (10-1)/3=3. A node with a multiple of three does not generate new $\theta \cdot 2^n$ sequences, so we then double the number 10: 10→20→ (40-1)/3=13. The node with the number 13 is not a multiple of three, so it generates a new $13\cdot 2^m$ sequence, with which we build a new branch on the *JT*, doubling the number 13: 13→26→ 52→ (52-1)/3=17.

The node with the number 17 is not a multiple of three, so it generates a new $17\cdot 2^m$ sequence, with which we build a new branch on the *JT*, doubling the number 17: 17→ (34-1)/3=11. The node with the number 11 is not a multiple of three, so it generates a new $11\cdot 2^m$ sequence, with which we build a new branch on the *JT*, doubling the number 11: 11→ 22→ (22-1)/3=7. The node with the number 7 is not a multiple of three, so it generates a new $7\cdot 2^m$ sequence, with which we build a new branch on the *JT*, doubling the number 7: 7→ 14→ (28-1)/3=9. A node with the number 9 is a multiple of three does not generate new $\theta \cdot 2^n$ sequences, so we then double the number 28: 28→ 56→ (112-1)/3=37 and i.e. On the $13 \cdot 2^m$ branch of the Jacobsthal number 17 form nodes: 17,69,277, ... Therefore, by similarly doubling the number 52, from node 17 we build a new branch: 52→ 104→ 208 (208-1)/3=69. A node with the number 99 is a multiple of three does not generate new $\theta \cdot 2^n$ sequences, so we then double the number 208: 208→ 416→ 832 → (832-1)/3=277 and i.e. So, from an arbitrary node with a non-multiple three number, all other branches of the *JT* are built. The algorithm for constructing the *JT* for the Collatz function $C_{3q-1}$ is implemented similarly.

**Appendix 3.** *Illustration of the correlation nods* $m(p)_{3,\theta_{1,5},n}$ *of JT with* $CT_{3q\pm1}$

Using formulas

$$\theta_{1,5} \cdot 2^{r(s)} = 3m_{\theta_{1,5},r(s)} \pm 1 \qquad \text{(A3-1)}$$

we construct the *JT*, which we reconstruct $CT_{3q\pm1}$ in the reverse direction $n \to 0$. To do this, we will create cells of the type

$$\begin{matrix} & p_{\theta_1,1} & & p_{\theta_1,3} & & p_{\theta_1,5} & & \\ (\theta_1 2^0 & \theta_1 2^1 & \theta_1 2^2 & \theta_1 2^3 & \theta_1 2^4 & \theta_1 2^5 & \theta_1 2^6 & \ldots) \\ m_{\theta_1,0} & & m_{\theta_1,2} & & m_{\theta_1,4} & & m_{\theta_1,6} & \end{matrix} \Rightarrow \begin{cases} \theta \cdot 2^{N_{even}} = m_{3,\theta,N_{even}} + p_{3,\theta,N_{even}+1}, \\ \theta \cdot 2^{N_{odd}} = p_{3,\theta,N_{odd}} + m_{3,\theta,N_{odd+1}}, \end{cases} \qquad \text{(A3-2)}$$

and

$$\begin{matrix} p_{\theta_5,0} & & p_{\theta_5,2} & & p_{\theta_5,4} & & p_{\theta_5,6} \\ (\theta_5 2^0 & \theta_5 2^1 & \theta_5 2^2 & \theta_5 2^3 & \theta_5 2^4 & \theta_5 2^5 & \theta_5 2^6 \ldots) \\ & m_{\theta_5,1} & & m_{\theta_5,3} & & m_{\theta_5,5} & \end{matrix} \Rightarrow \begin{cases} \theta \cdot 2^{N_{even}} = p_{3,\theta,N_{even}} + m_{3,\theta,N_{even}+1}, \\ \theta \cdot 2^{N_{odd}} = m_{3,\theta,N_{odd}} + p_{3,\theta,N_{odd+1}}, \end{cases} \text{(A3-3)}$$

where the numbers $m(p)_{\theta_{1,5},N_{even(odd)}}$ in cells (A3-2) and (A3-3) form nodes. The nodes of $\theta_3 \cdot 2^n$ are active and of type $\theta_3 \cdot 2^n$ inactive. The members of inactive cells $\theta_3 \cdot 2^n$ are formed from duplicated values $\theta_3 \cdot 2^n$.

$$\begin{array}{l} \to \overset{181}{\downarrow} \to \qquad \theta_1 = 181 \\ \to 544 \to 438 \to 136 \to 68 \to 34 \to \overset{17}{\downarrow} \to \qquad \theta_5 = 17 \\ \qquad \to 52 \to 26 \to \overset{13}{\downarrow} \to \qquad \theta_1 = 13 \\ \qquad\qquad \to 40 \to 20 \to 10 \to \overset{5}{\downarrow} \to \qquad \theta_5 = 5 \\ \qquad\qquad\qquad \to 16 \to 8 \to 4 \to 2 \to 1 \qquad \theta_{END} = 1 \end{array}$$

**Figure A3-1.** *Graph conjecture* $C_{3q+1}$ *of the number* $q = 181$

However, as can be seen from (A3-1), the Jacobsthal numbers $K^{\pm}_{\theta_{1,5},n}$ can be structured in the form of cells in which the active nodes $\theta_1 = m(p)_{\theta_1,r(s)}$ and $\theta_5 = m(p)_{\theta_5,r(s)}$ in the direction of increasing numbers $m(p)_{\theta_{1,5},r(s)}$ from $\theta_1$ to $\theta_5$ are surrounded by two multiples of three numbers $[m(p)_{\theta_3,r(s)}]$. For example, the first cell of the polynomial $5\cdot 2^n$ looks like this:

$$\begin{array}{cccc} & \text{type } \theta_1 & \text{type } \theta_5 & \\ & \Downarrow & \Downarrow & \\ \text{—}[3]\text{—} & 13 & \text{—}\ 53\ \text{—} & [213\ ]\text{—} \\ \Uparrow & \text{I} & \text{I} & \Uparrow \\ \text{type } \theta_3 & & & \text{type } \theta_3 \end{array} \qquad \text{(A3-4)}$$

On the sequ $\cdot 2^n$ the cells ( is equal:

$$\begin{array}{cccccccc} & p_{1,1} & & p_{1,3} & & p_{1,5} & & \\ (1\cdot 2^0 & 1\cdot 2^1 & 1\cdot 2^2 & 1\cdot 2^3 & 1\cdot 2^4 & 1\cdot 2^5 & 1\cdot 2^6 & ...) \\ m_{1,0} & & m_{1,2} & & m_{1,4} & & m_{1,6} & \end{array} \qquad \text{(A3-5)}$$

Therefore. the first cell is equal

$$\theta_1 : \left\{ \begin{array}{ccccccc} m_{1,0} & & m_{1,2} & & m_{1,4} & & m_{1,6} \\ [0] & - & 1 & - & 5 & - & [21]- \end{array} \right. \qquad \text{(A3-6)}$$

is formed by the first four of the Jacobsthal numbers $0,\ 1,\ 5,\ 21,\ 85,\ 341,\ 1024,\ ...$ and has only one active node $m_{3,1,4} = 5$, which $CT_{3q+1}$ leads to the unit value (Figure A3-1).

For the sign-symmetric transformation $C_{3q-1}$, the first cell is equal

$$\theta_1 : \left\{ \begin{array}{ccccccc} p_{1,1} & & p_{1,3} & & p_{1,3} & & p_{1,5} \\ [3] & - & 11 & - & 43 & - & [171]- \end{array} \right. \qquad \text{(A3-7)}$$

is formed by the first four of the numbers : 1, 3, 11, 43, 171, 683, 2731, 10923,...and has two active node $p_{1,5} = 11$ and $p_{1,7} = 43$, which $CT_{3q-1}$ leads to the unit value (Figure A3-2).

$$\rightarrow \overset{79}{\underset{\downarrow}{}} \rightarrow$$

$$\rightarrow 236 \rightarrow 118 \rightarrow \overset{59}{\underset{\downarrow}{}} \rightarrow$$

$$\rightarrow 176 \rightarrow 88 \rightarrow 44 \rightarrow 22 \rightarrow \overset{11}{\underset{\downarrow}{}} \rightarrow$$

$$\rightarrow 32 \rightarrow 16 \rightarrow 8 \rightarrow 4 \rightarrow 2 \rightarrow 1 \qquad \theta_{END} = 1$$

**Figure A3-2.** *Graph conjecture* $C_{3q-1}$ *of the number* $q = 79$

**Appendix 4.** *Example of calculating numbers* $K^{\pm}_{3,\theta_1,n}$ *by formula (1.5) from* $\theta = 3\cdot I$

| n | 0 | 1 | 2 | 3 | 4 | 5 | … | |
|---|---|---|---|---|---|---|---|---|
| $\theta$=3 | 4/3 | 5/3 | 13/3 | 23/3 | 49/3 | 35/3 | … | $K^{-}_{3,\theta,n}$ |
| $\theta$=3 | 2/3 | 8/3 | 11/3 | 25/3 | 47/3 | 97/3 | … | $K^{+}_{3,\theta,n}$ |
| $\theta$=9 | 10/3 | 17/3 | 37/3 | 71/3 | 145/3 | 287/3 | … | $K^{-}_{3,\theta,n}$ |
| $\theta$=9 | 8/3 | 19/3 | 35/3 | 73/3 | 143/3 | 289/3 | … | $K^{+}_{3,\theta,n}$ |
| $\theta$=15 | 16/3 | 29/3 | 61/3 | 119/3 | 241/3 | 479/3 | … | $K^{-}_{3,\theta,n}$ |
| $\theta$=15 | 14/3 | 31/3 | 59/3 | 121/3 | 239/3 | 481/3 | … | $K^{+}_{3,\theta,n}$ |

(A4-1)

Thus, if the $\theta$ index of the $\theta \cdot 2^n$ sequence is a multiple of three, the $K^{\pm}_{3,\theta_1,n}$ numbers take fractional values. The numbers $K^{\pm}_{3,\theta_1,n}$ are fractal if $\theta = a$, unless $a = 1$.

**Appendix 5.** Let's convert numbers from the interval of $q = 1 \div 300{,}000$ values with the $C_{3q-1}$ function and divide them by $cycle^{3q-1}_{1\leftrightarrow 1}, cycle^{3q-1}_{5\leftrightarrow 5}$ and $cycle^{3q-1}_{17\leftrightarrow 17}$ cycles. We find that of the frequencies $\nu = \dfrac{Number_{cycle}}{Number_{full}}$

$$\nu_{cycle^{3q-1}_{1\leftrightarrow 1}} \cong 0.33, \quad \nu_{cycle^{3q-1}_{5,7\leftrightarrow 5,7}} \cong 0.3{,}2 \quad \nu_{cycle^{3q-1}_{17\leftrightarrow 17}} \cong 0.35 \qquad \text{(A5-1)}$$

distribution on the cycles of $cycle^{3q-1}_{1\leftrightarrow 1}$, $cycle^{3q-1}_{5\leftrightarrow 5}$ and $cycle^{3q-1}_{5\leftrightarrow 5}$ :

$$cycle^{3q-1}_{1\leftrightarrow 1} = \{1 \leftrightarrow 2 \leftrightarrow 1 \leftrightarrow 2\}, \quad cycle^{3q-1}_{5\leftrightarrow 5} = \{5 \leftrightarrow 14 \leftrightarrow 7 \leftrightarrow 20 \leftrightarrow 10 \leftrightarrow 5\} \qquad \text{(A5-2)}$$

$$cycle^{3q-1}_{17\leftrightarrow 17} = \{17 \leftrightarrow 50 \leftrightarrow 25 \leftrightarrow 74 \leftrightarrow 37 \leftrightarrow 110 \leftrightarrow 55 \leftrightarrow 164 \leftrightarrow 82 \leftrightarrow 411 \leftrightarrow$$
$$\leftrightarrow 122 \leftrightarrow 61 \leftrightarrow 182 \leftrightarrow 91 \leftrightarrow 272 \leftrightarrow 136 \leftrightarrow 68 \leftrightarrow 34 \leftrightarrow 17\}$$

**Table A5-1.** *Illustration a fragment of the distribution of the set of the first 200 natural numbers by three isolated cycles* $cycle^{3q-1}_{1\leftrightarrow 1}$, $cycle^{3q-1}_{5\leftrightarrow 5}$, $cycle^{3q-1}_{17\leftrightarrow 17}$ *from the 3q-1 conjecture*

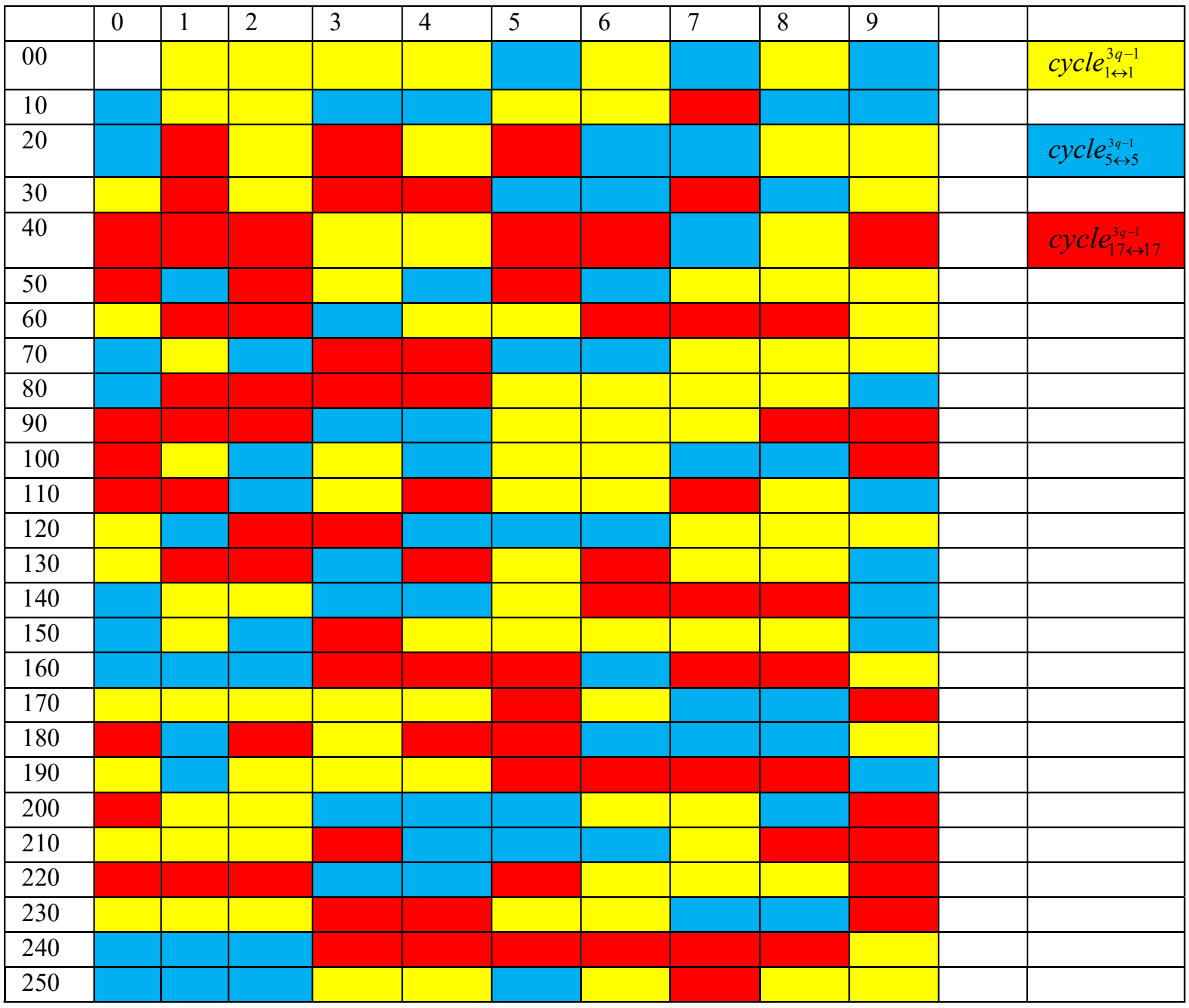

Therefore, from conjecture $3q-1$ the frequency $\nu$ the trend toward a uniform distribution. Then from the point of view of the Collatz process, the cycles $cycle^{3q-1}_{1\leftrightarrow 1}$, $cycle^{3q-1}_{5(7)\leftrightarrow 5(7)}$, $cycle^{3q-1}_{17\leftrightarrow 17}$ are equivalent to each other.

The $CT_{3q-1}$ of numbers $q_5 = 5, q_7 = 7, q_{17} = 17$,do not reach a single value $1 \cdot 2^0$, but have the form of final periodic oscillations with minimal amplitudes $q_{\min} = 5, 7, 17$. A triple of recurrent numbers $q_5, q_7, q_{17}$ satisfies the additional condition:

$$q_5 \cdot q_7 = 2q_{17} + 1 \Rightarrow J_k \cdot LJ_{k+1} = 2LJ_{k+2} + 1, \qquad \text{(A5-3)}$$

where its first element is equal to:

$$(2^k + 1)(2^{k+1} - 1) = 2^{k+2} + 1 \Rightarrow 2^k = 4 \Rightarrow k = 2. \qquad \text{(A5-4)}$$

In the direction $k \to \infty$, the product grows exponentially, so there are no other triples of numbers with properties (A5-3). The triplet of numbers $q_5, q_7, q_{17}$ with properties (A5-4) are absent in $3q+1$ conjecture. Therefore, there are no isolated from the polynomial $1 \cdot 2^n$ in the conjecture $3q+1$.In addition, the numbers $q_5, q_{17}$ are prime Fermat numbers

$$\Phi_n = 2^{2^n} + 1 \Rightarrow \Phi_1 = 2^{2^1} + 1 = 5 \text{ and } \Phi_2 = 2^{2^2} + 1 = 17, \qquad \text{(A5-5)}$$

and the number $q_7$ is the Mersenne numbers $M_{n=3} = 2^3 - 1 = 7$. The $q_5, q_7, q_{17}$ numbers are also connected by the relation:

$$M_2 \Phi_3 - 2\Phi_4 = 1 \qquad \text{(A5-6)}$$

## Appendix 6

Figure A6-1a shows the so-called flat Jacobsthal lattice, the nodes of which are

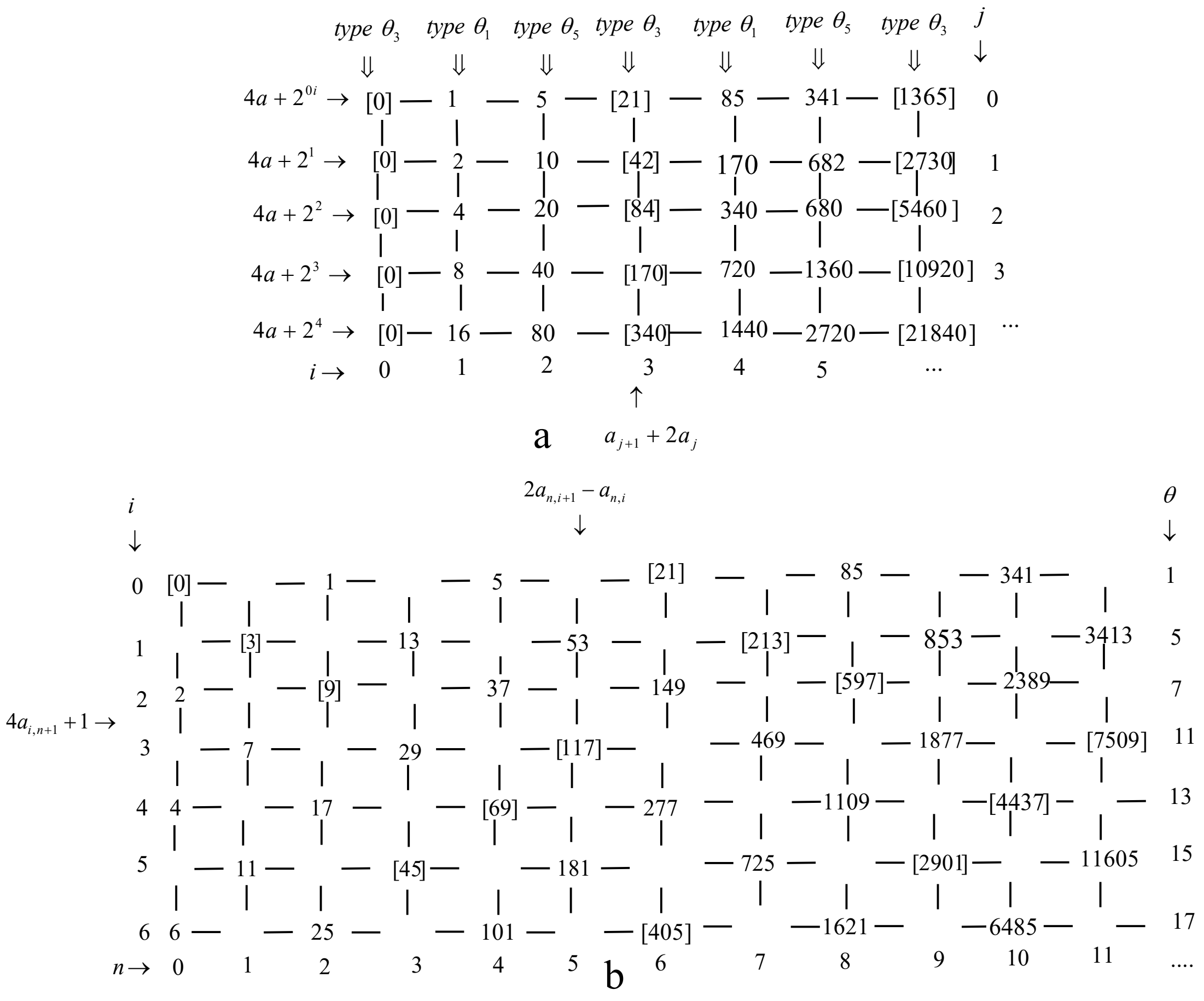

**Figure A6-1**. *Jacobsthal lattice from 3q+1 conjecture*

formed by transformations $4a_{i,j} \pm 2^j$ between adjacent members in rows and $a_{i,j+1} + 2a_{i,j}$ between adjacent members in columns. Analogous lattices hold true for Jacobsthal numbers $K^{\pm}_{\theta_5,n}$ with arbitrary indices $\theta$ (Figure A6-1b).

Figure A6-2 shows the so-called flat Collatz lattice, the nodes of which are formed by transformations $\theta \cdot 2^n$ between adjacent members in rows and $2^m(3q \pm 1)$ between adjacent members in columns. The lattice is semi-bounded, the right side of which is bounded $CT_{3q+1}$ (the trajectory is shown by arrows in the figure), and the opposite side is multiplied to infinity. The other two opposite sides are formed by initial $\theta \cdot 2^n$ and final $\Omega \cdot 2^m$ polynomials. The grid in Figure A6-2 shows that classical of $CT_{3q+1}$ are the optimal algorithm for achieving a single value by the Collatz sequence.

Figure A6-3 shows an illustration of the effect of completion in the reverse direction $n \to \infty$ of the process of filling Collatz sequences $C_{3q+1}$ with an odd number that is a multiple of three. These graphs are built on the basis of the rules for transformations of Jacobsthal numbers

$$m_{\theta_{1,5},r(s)} = \frac{\theta_{1,5} \cdot 2^{r(s)} - 1}{3} = I \text{ (a)} \quad \text{and} \quad p_{\theta_{1,5},r(s)} = \frac{\theta_{1,5} \cdot 2^{r(s)} + 1}{3} = I \text{ (b)} \qquad \text{( A6-1)}$$

As can be seen, sequences of odd Collatz numbers on the left are bounded by one, and on the right they lead to a multiple of three odd numbers $[odd_3]$. For example:

| | | | | | | | | | | | | | | |
|---|---|---|---|---|---|---|---|---|---|---|---|---|---|---|
| 3 | 5 | 1 | | | | | | | | | | | | |
| 9 | 7 | 11 | 17 | 13 | 5 | | | | | | | | | |
| 15 | 23 | 35 | 53 | 5 | | | | | | | | | | |
| 21 | 1 | | | | | | | | | | | | | |
| 27 | 41 | 31 | 47 | 71 | 107 | 161 | 121 | 91 | 137 | 103 | 155 | 233 | 175 | 263 |
| | 395 | 593 | 445 | 167 | 251 | 377 | 283 | 425 | 319 | 479 | 719 | 1079 | 1619 | 2429 |
| | 911 | 1367 | 2051 | 3077 | 577 | 433 | 325 | 61 | 23 | | | | | |
| 33 | 25 | 19 | 29 | 11 | | | | | | | | | | |
| 39 | 59 | 89 | 67 | 101 | 19 | | | | | | | | | |
| 45 | 17 | | | | | | | | | | | | | |
| 51 | 77 | 29 | | | | | | | | | | | | |
| 57 | 43 | 65 | 49 | 37 | 7 | | | | | | | | | |
| 63 | 95 | 143 | 215 | 323 | 485 | 91 | | | | | | | | |
| 69 | 13 | | | | | | | | | | | | | |
| 75 | 113 | 85 | 1 | | | | | | | | | | | |
| 81 | 61 | 23 | | | | | | | | | | | | |
| 87 | 131 | 197 | 37 | | | | | | | | | | | |
| 93 | 35 | | | | | | | | | | | | | |
| 99 | 149 | 7 | | | | | | | | | | | | |
| 105 | 79 | 119 | 179 | 269 | 101 | 19 | | | | | | | | |
| 111 | 167 | | | | | | | | | | | | | |
| 117 | 11 | | | | | | | | | | | | | |
| 123 | 185 | 139 | 209 | 157 | 59 | | | | | | | | | |
| 129 | 97 | 73 | 55 | 83 | 125 | 47 | | | | | | | | |
| 135 | 203 | 305 | 229 | 43 | 65 | | | | | | | | | |
| 141 | 53 | | | | | | | | | | | | | |
| 147 | 221 | 83 | 125 | 47 | | | | | | | | | | |

or in the general case:

$\{1, ..., [odd_3]\}$. (A6-2)

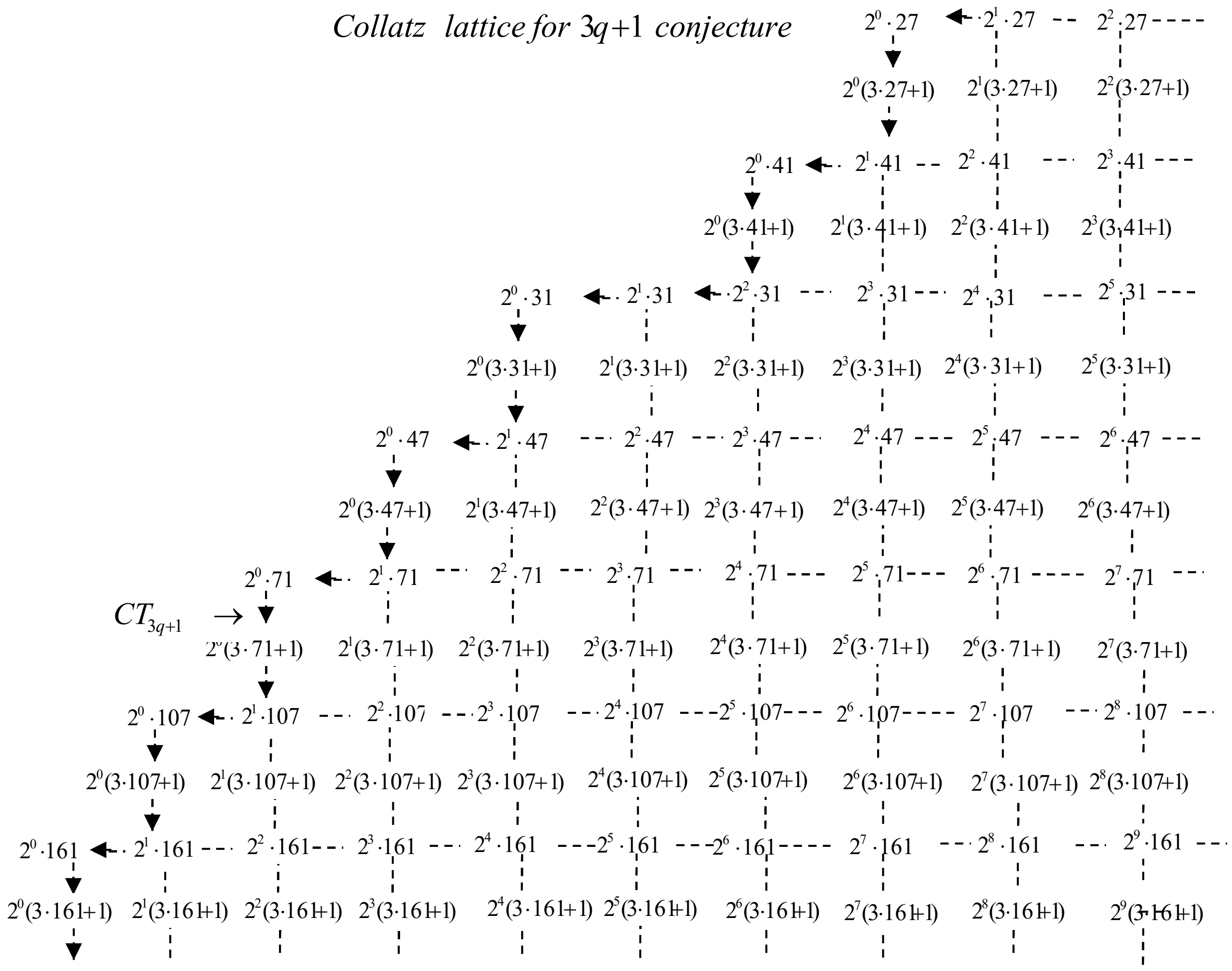

**Figure A6-2.** *Collatz lattice for 3q+1 conjecture*

To conjecture (2), an arbitrary full graph $Graph_{full}$, can be considered as a superposition of two types of graphs:

$$Graph_{full} = Graph\ I + Graph\ II, \tag{A6-3}$$

where the graph $Graph\ I$ is formed by a sequence of odd numbers of intermediate calculations $X_n$:

$$Graph\ I = \left\{1, ..., \frac{a_N - 1}{3} = Od_3\right\} \tag{A6-4}$$

the left boundary of which is an infinite loop infinite cycle $cycle^{3n+1}_{1\to 4\to 1} = \{1 \to 4 \to 2 \to 1\}$.

The right boundary (A6-4) is formed by an odd number that is a multiple of three $X_n = Od_3$. On graph $Graph\ I$ nodes are formed, from which side graphs are subsequently generated, so they can be considered as active. The graph $Graph\ II$, generated by a node with a number $\frac{a_N - 1}{3} = Od_3$, and then in the reverse direction is formed from double values $X_n = Od_3$:

$$Graph\ II = 2^0 \cdot Od_3, \quad 2^1 \cdot Od_3, \quad 2^2 \cdot Od_3, \quad 2^3 \cdot Od_3, \ldots \qquad \text{(A6-5)}$$

*Collatz conjecture of odd numbers from* 3 *to* 189

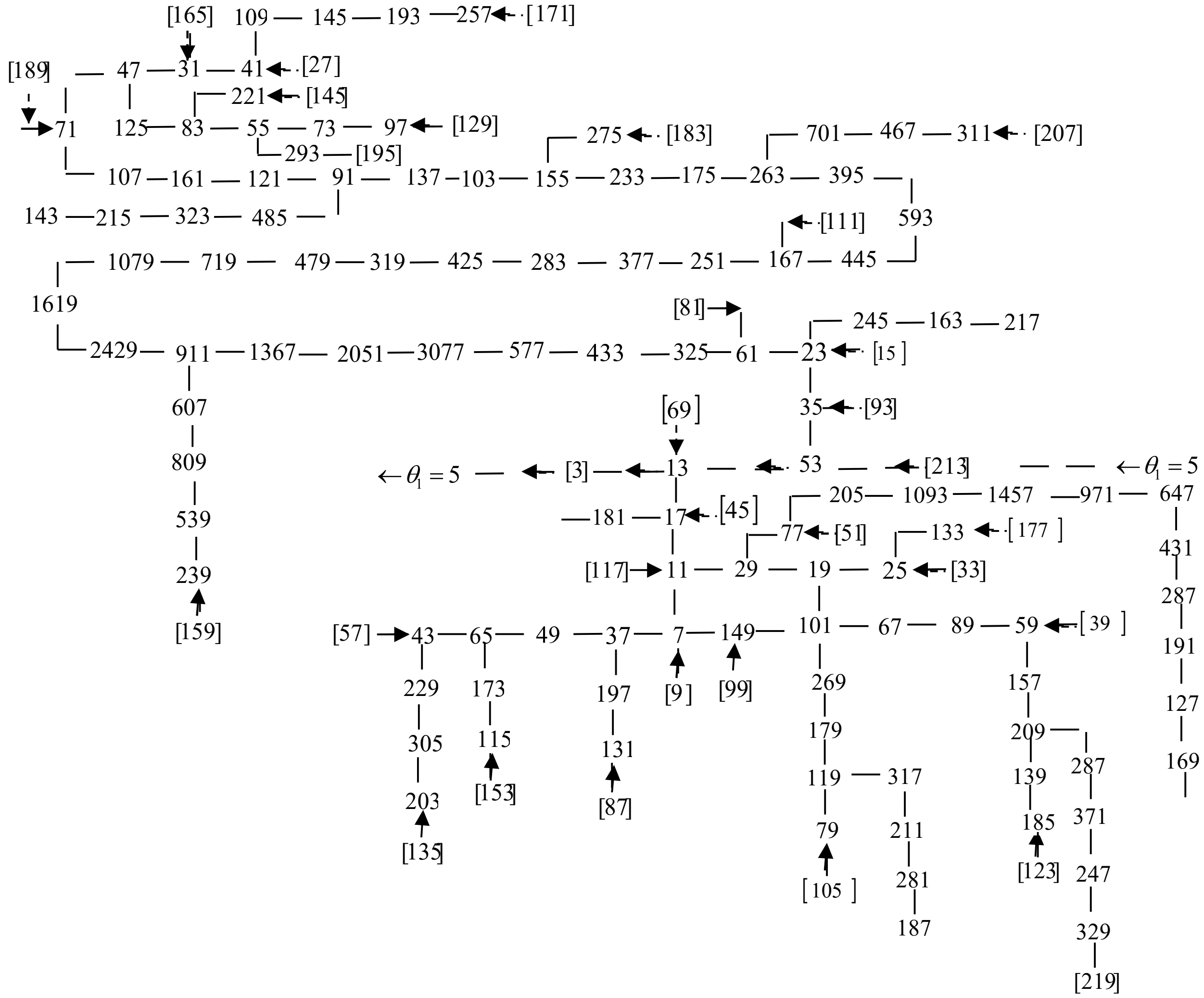

**Figure A6-3**. *Collatz conjecture of odd number from 3 to 189*

In each sequence, the odd numbers that are not a multiple of three are the sources of the formation of new Collatz sequences, each of which ends with an odd number that is a multiple of three. On the *JT* , all odd numbers are sources of formation of other sequences of type $\theta \cdot 2^m$. Sequences $\theta \cdot 2^m$ are of two types, without branches if $\theta = odd_3$ is a multiple of three, and with branches for numbers $\theta \neq odd_3$. Moreover, in each interval of the type from 1 to 199 out of 100 odd numbers, there are 32 numbers of the type $\theta = odd_3$ and 68 numbers of the type $\theta \neq odd_3$.

**Acknowledgements**

The author thanks B. M Gurbaxani and D.Barina for interesting and useful discussions. The author expresses his gratitude to the members of the seminar of the Department of Higher Mathematics of the National University of Lviv Polytechnic and The seminar of the Chair of Algebra, Topology and Foundations of Mathematics of the Franko Lviv University, and seminar on fractal analysis (online) of the Institute of Mathematics of NAS of Ukraine, Kyiv, and V. Ilkiv, T. Banakh, M. Samotiy, M. Pratsiovytyi, D.Karvatskyi, O. Baranovskyi for the advice and critical remarks expressed during the discussion of the article. The anonymous referee's detailed comments and suggestions were also very valuable for references and the quality of exposition of this paper.

**Conflicts of Interest**

The author declares no conflicts of interest regarding the publication of this paper

**Funding**

This research was carried out without any funding